\documentclass{amsart}
\usepackage{amssymb,latexsym}
\usepackage{amsfonts}

\newcommand{\ba}{\begin{array}}
\newcommand{\ea}{\end{array}}

\def \qed{\cqfd}

\def\qed{\vbox{\hrule
\hbox{\vrule\hbox to 5pt{\vbox to 8pt{\vfil}\hfil}\vrule}\hrule}}

\newcommand{\beg}{\begin{eqnarray*}}
\newcommand{\begn}{\begin{eqnarray}}
\newcommand{\en}{\end{eqnarray*}}
\newcommand{\enn}{\end{eqnarray}}

\begin{document}
\title{Energy properness and Sasakian-Einstein metrics}
\author{Xi Zhang}
\address{Department of Mathematics\\
Zhejiang University, P. R. China } \email{ xizhang@zju.edu.cn}
\thanks{The work was supported in part by NSF in
China, No.10831008.}

\begin{abstract}
In this paper,  we  show that the existence  of
Sasakian-Einstein metrics is closely related to the properness of corresponding energy functionals. Under the condition that admitting no nontrivial Hamiltonian holomorphic vector field, we prove that the existence of Sasakian-Einstein metric implies a Moser-Trudinger type inequality. At the end of this paper, we also obtain a Miyaoka-Yau type inequality in Sasakian geometry.

\end{abstract}

\maketitle

\section{Introduction}
\setcounter{equation}{0}

An odd dimensional Riemannian manifold $(M , g)$ is said to be a
Sasakian manifold if the cone manifold $(C(M) , \tilde{g} ) =(
M\times R^{+} , \quad r^{2}g + dr^{2})$ is K\"ahler. In this paper, we suppose that $dim M =2m+1$. Furthermore, Sasakian manifold $(M, g)$ is said to be
Sasakian-Einstein if the Ricci tensor of $g$
 satisfies the Einstein condition. It is well
known that the K\"ahler cone $(C(M) , \tilde{g} )$ must be a Calabi-Yau cone if $(M,
g)$ is a Sasakian-Einstein manifold. Recently, Sasakian-Einstein metrics have attract increasing attention, as they provide rich
source of constructing new Einstein manifolds in odd dimensions
 and its important role in the superstring theory, see references \cite{1, 2, 3, Boyer, 4, 5, BGS, GMSW1, GMSW2, GMSY, GZ, 8, 12, 15, 16, 17, 18, Z1}.

 Sasakian manifolds can be
studied from many view points as they have many structures. A
Sasakian manifold $(M, g)$ has a contact structure $(\xi , \eta ,
\Phi )$, and it also has a one dimensional foliation $\mathcal
F_{\xi}$, called the Reeb foliation.  Here,  the Killing vector field $\xi $ is
called the characteristic or Reeb vector field, $\eta $ is called
the contact $1$-form, $\Phi$ is a $(1 , 1)$ tensor field which
defines a complex structure on the contact sub-bundle $\mathcal
D=ker \eta $. In the following, a Sasakian manifold will be denoted
by $(M, \xi , \eta , \Phi , g)$, the quadruple $(\xi , \eta , \Phi
, g)$ will be called by a Sasakian structure on manifold $M$. In a natural way, a Sasakian structure
$(\xi , \eta , \Phi
, g)$ induce a transverse holomorphic structure and a transverse K\"ahler metric on the foliation $\mathcal
F_{\xi}$. In this paper, we may change Sasakian structure, but always fix the Reeb vector field $\xi $ and the transverse
holomorphic structure on $\mathcal F_{\xi}$.

   Fixed a transverse holomorphic structure on $\mathcal F_{\xi}$,  we have a splitting of the complexification of the
bundle $\wedge ^{1}_{B}(M)$ of basic one forms on $M$, $ \wedge
^{1}_{B}(M) \otimes C =\wedge ^{1, 0}_{B}(M)\oplus \wedge ^{0,
1}_{B}(M)$, and then we  have the decomposition of $d$, i.e.
$d=\partial _{B}+\bar{\partial }_{B}$. We  have the basic
cohomology groups $H^{i, j}_{B}(M, \mathcal F_{\xi})$ which enjoy
many of the same properties as the Dolbeault cohomology of a
K\"ahler structure. We also have the transverse Chern-Weil theory and can define the basic Chern classes $c_{k}^{B}(M, \mathcal F_{\xi })$. For the detail, see \cite{Boyer}. Given a transverse K\"ahler structure $g^{T}$,
 one can define the transverse
Levi-Civita connection $\nabla^{T}$ on the normal bundle $\nu (\mathcal F _{\xi }) =TM/L\xi$, and then one
can define the transverse Ricci curvature $Ric^{T}$, see section 2
for details. We denote the related Ricci form by $\rho^{T}$, it is
easy to see that $\rho^{T}$ is a closed basic $(1, 1)$-form and
the basic cohomology class $\frac{1}{2\pi}[\rho
^{T}]_{B}=c_{1}^{B}(M , \mathcal F _{\xi })$ is the basic first
Chern class.
A Sasakian metric $(\xi , \eta , \Phi , g)$ is said
to be transversely K\"ahler-Einstein if it's transverse Ricci form satisfies
$\rho^{T}=\mu d\eta $.  It's easy to see that a Sasakian metric
$(\xi , \eta , \Phi , g)$ is Sasakian-Einstein then it must be
transversely K\"ahler-Einstein and $\rho^{T}=(m+1) d\eta $. So, a necessary
condition for the existence of Sasakian-Einstein metric on $M$ is
that there exists a Sasakian structure $(\xi , \eta , \Phi , g)$
such that $2\pi c_{1}^{B}(M, \mathcal F_{\xi })=(m+1)[d\eta]_{B}$.

\medskip

Given a Sasakian structure $(\xi , \eta , \Phi , g)$ on $M$,  Let's denote the space of all smooth  basic real function $\varphi
$ (i.e. $\xi \varphi \equiv 0$) on $(M, \xi , \eta , \Phi , g)$ by
$C_{B}^{\infty}(M, \xi)$. Set
\begin{eqnarray}\label{H}
\mathcal H (\xi , \eta , \Phi , g)=\{\varphi \in C_{B}^{\infty}
(M, \xi) : \eta_{\varphi }\wedge (d \eta_{\varphi })^{n} \neq 0
\},
\end{eqnarray}
where
\begin{eqnarray}\label{eta-d}
\eta_{\varphi } =\eta +\sqrt{-1}\frac{1}{2}(\bar{\partial
}_{B}-\partial_{B})\varphi , \quad d\eta_{\varphi } =d\eta +
\sqrt{-1 }\partial_{B} \bar{\partial }_{B}\varphi .
\end{eqnarray}
 For any $\varphi \in
\mathcal H$, $(\xi , \eta _{\varphi } , \Phi_{\varphi },
g_{\varphi }) $ is also a Sasakian structure on $M$, where
\begin{eqnarray}\label{phi-d} \Phi_{\varphi } = \Phi -\xi \otimes
(d_{B}^{c}\varphi ) \circ \Phi ,\quad g_{\varphi } =\frac{1}{2}
d\eta_{\varphi } \circ (Id \otimes \Phi_{\varphi } ) +
\eta_{\varphi } \otimes \eta_{\varphi }.
\end{eqnarray}
Furthermore,  $(\xi , \eta _{\varphi } , \Phi_{\varphi }, g_{\varphi }) $ and
$(\xi , \eta , \Phi , g)$ have the same transversely holomorphic
structure on $\nu (\mathcal F _{\xi})$ and the same holomorphic
structure on the cone $C(M)$ (Proposition 4.2 in \cite{12}, also
\cite{Boyer} ). Obviously, those deformations of Sasakian
structure deform the transverse K\"ahler form in the same basic
$(1, 1)$ class. We call this class the basic K\"ahler class of the
Sasakian manifold $(M, \xi , \eta , \Phi , g)$.

As in the K\"ahler case, one can define Aubin's functionals $I_{d\eta }$, $J_{d\eta }$,   Ding and Tian's energy functional $F_{d\eta }$, and Mabuchi's $K$-energy functional $\mathcal V_{d\eta }$  on the space $\mathcal H (\xi , \eta , \Phi , g)$, see section 3 for details. We say the energy functional $F_{d\eta }$ (or $\mathcal V_{d\eta }$) is proper if $\limsup _{i\rightarrow +\infty} F_{d\eta } (\varphi_{i})=+\infty$ whenever
$\lim _{i\rightarrow +\infty} J_{d\eta } (\varphi_{i})=+\infty$, where $\varphi_{i}\in \mathcal H (\xi , \eta , \Phi , g)$.
We say two Sasakian structures are compatible with each other if they have the same Reeb vector field and the same transverse holomorphic structure.

In K\"ahler geometry, Tian \cite{T3} have show that  the existence of K\"ahler-Einstein metric is equivalent with
the properness of $F$ energy functional on a compact K\"ahler manifold with positive Chern class and without any nontrivial holomorphic fields.
 In this paper, under the condition that without nontrivial Hamiltonian holomorphic vector field (see section 2 for details), we generalize Tian's result in \cite{T3} to the Sasakian case. In fact, we obtain the following theorem.

\medskip

{\bf Main Theorem } {\it Let  $(M, \xi , \eta , \Phi , g)$ be a compact Sasakian manifold with $[d\eta ]_{B} = \frac{2\pi }{m+1}
c_{1}^{B}(M, \mathcal F_{\xi})$ and without any nontrivial Hamiltonian holomorphic vector field. Then $M$ has a Sasakian-Einstein structure compatible with $(\xi , \eta , \Phi , g)$ if and only if the functional $F_{d\eta }$ (or the $K$-energy functional $\mathcal V_{d\eta }$) is proper in the space $\mathcal H(\xi , \eta , \Phi , g)$.}

\medskip

We will follow Tian's method, and the discussion in \cite{PSSW} by Phong, Song, Strum and Weinkove. It's easy to see that if a Sasakian structure
$(\xi ' , \eta ' , \Phi ', g ' ) $ is compatible with
$(\xi  , \eta  , \Phi  , g )$, then we have $[d\eta ']_{B}=[d \eta ]_{B}\in H_{B}^{1, 1}(\mathcal F_{\xi })$, by transverse $\partial \bar{\partial }$ lemma (\cite{11}), there exists a basic function $\varphi \in \mathcal H (\xi , \eta , \Phi , g)$ such that
\begin{eqnarray}
d\eta '=d\eta +dd^{c}_{B}\varphi , \quad and \quad \eta '=\eta +d_{B}^{c}\varphi +\zeta
\end{eqnarray}
where $\zeta $ is closed basic one form, $d^{c}_{B}=\frac{\sqrt{-1}}{2}(\bar{\partial }_{B}-\partial_{B})$. So, the existence problem of Sasakian-Einstein metric compatible with $(\xi , \eta , \Phi , g)$ can be reduced
to solving the following transverse Monge-Amp\`ere equation,
\begin{eqnarray}\label{MA1}
\frac{(d\eta+\sqrt{-1}\partial_{B}\overline{\partial }_{B}\varphi
)^{m}\wedge \eta}{(d\eta )^{m}\wedge \eta}=\exp (h_{d\eta}-(m+1)\varphi ),
\end{eqnarray}
where $\varphi \in \mathcal H (\xi  , \eta  , \Phi  , g )$ and $h_{d\eta}$ is a smooth basic function
which satisfies $\rho^{T}=(m+1)d\eta
+\sqrt{-1}\partial_{B}\bar{\partial }_{B}h_{d\eta}$ and $\int_{M}\exp(h_{d\eta}) (d\eta )^{m}\wedge \eta =\int_{M} (d\eta )^{m}\wedge \eta=V$. In order to use the
continuity method, we consider the following family of equations
\begin{eqnarray}\label{MA2}
\frac{(d\eta+\sqrt{-1}\partial_{B}\overline{\partial }_{B}\varphi
)^{m}\wedge \eta}{(d\eta )^{m}\wedge \eta}=\exp (h_{d\eta}-t(m+1)\varphi
),
\end{eqnarray}
 where  $t\in [0,
 1]$. By El-Kacimi's ( \cite{11})
generalization of Yau's estimate (\cite{Y}) for transverse
Monge-Amp\`ere equations, to solve the transverse Monge-Amp\`ere
equation (\ref{MA1}), it is sufficient to obtain a priori uniform $C^{0}$
estimate of solution $\varphi $ of the transverse Monge-Amp\`ere
equation (\ref{MA2}). Different
than the K\"ahler case, the transverse Monge-Amp\`ere equation
(\ref{MA2}) only give bounds on the transverse Ricci curvature
which does not lead a lower bound of the Ricci curvature by a
positive constant. So, we can not apply the Myers theorem directly
to obtain an estimate on the diameter and the lower bound of the Green's function. It's pointed out by Sekiya (\cite{Se}),
through $\mathcal D$-homothetic deformation of Sasakian structure,
one can get the desired estimates in the Sasakian case. In section 4 (theorem 4.3), we show that the properness of energy functional $F_{d\eta}$ (or the $K$-energy functional $\mathcal V_{d\eta }$) implies the $C^{0}$ estimates, and so we get a existence result for Sasakian-Einstein metric.

 To prove the existence of Sasakian-Einstein metric implies the properness of energy functional, we use the backward continuity method. In order to apply the implicity function theorem at time $t=1$, we should consider the eigenspace  corresponding to the eigenvalue $-4(m+1)$ of the basic Laplacian of the Sasakian-Einstein metric. Thanks to  A. Futaki, H. Ono and G. Wang's result \cite{12}, we know that the above eigenspace is not empty if and only if there exists a nontrivial Hamiltonian holomorphic vector field.  In section 6 (theorem 6.1), under the condition that without nontrivial Hamiltonian holomorphic vector field, we obtain the following Moser-Trudinger inequality on Sasakian-Einstein manifold $(M, \xi , \eta_{SE}, \Phi_{SE}, g_{SE})$, i.e.
 there exist uniform positive constants $C_{1}$, $C_{2}$, such that
\begin{eqnarray}\label{1MT2}
F_{d\eta_{SE} }(\varphi) \geq C_{1}J_{d\eta_{SE} }(\varphi ) -C_{2},
\end{eqnarray}
for all $\varphi \in \mathcal H (\xi , \eta_{SE} , \Phi_{SE} , g_{SE})$. In view of the cocycle identity of $F_{d\eta }$ and properties of $J_{d\eta}$ (see section 3, lemma 3.1 and lemma 3.2), the inequality holds for every Sasakian structure $(\xi , \eta , \Phi , g)$ which compatible with the Sasakian-Einstein structure $(\xi , \eta_{SE} , \Phi_{SE} , g_{SE})$. On the other hand, the relation (\ref{MF3}) implies that the Moser-Trudinger type inequality (\ref{1MT2}) also be valid for the $\mathcal K$-energy $\mathcal V_{d\eta }$.

Given a Sasakian structure $(\xi, \eta, \Phi, g)$ on $M$ with $2\pi
c_{1}^{B}(M)=(m+1)[d\eta ]_{B}$, we denote $\mathcal S (\xi , \bar{J})$ to be the set of all Sasakian structures which compatible with $(\xi, \eta, \Phi, g)$.
By definition 2.6 and proposition 2.7, it's easy to see that the integral
\begin{eqnarray}\label{II}
\int_{M} (2c_{2}^{B}(M, \mathcal F_{\xi })-\frac{m}{m+1}c_{1}^{B}(M, \mathcal F_{\xi })^{2})\wedge \frac{(\frac{1}{2}d\eta')^{m-2}}{(m-2)!}\wedge \eta'
\end{eqnarray}
is independent of the choice of a Sasakian structure  $(\xi , \eta', \Phi', g')\in \mathcal S (\xi , \bar{J})$. By direct calculation, see lemma 7.2, if there exists a Sasakian-Einstein structure ( or equivalently a Sasakian structure with constant scalar curvature)in $\mathcal S (\xi , \bar{J})$, then the integral (\ref{II}) must be nonnegative, this inequality will be called a Miyaoka-Yau type inequality. In section 7, we discuss this Miyaoka-Yau type inequality under a more weak condition, and show that the energy function $F_{d\eta }$ (or Mabuchi's $\mathcal K$-energy $\mathcal V_{d\eta }$) bounded below will implies this weak condition, see theorem 7.3 and proposition 7.5 for details.

\medskip

This paper is organized as follows. In Section 2, we will recall
 some preliminary results about Sasakian geometry. In section 3,
 we introduce energy functionals in Sasakian
 geometry. In section 4, we consider the
 transverse Monge-Amp\`ere equation, and give a existence result of Sasakian-Einstein metric. In section 5, we use the Sasakian-Ricci flow to get a smoothing lemma. In  section 6, we give a proof of the Moser-Trudinger type inequality (\ref{1MT2}), and finish the proof of the main theorem. In the last section, we obtain a Miyaoka-Yau type inequality in Sasakian geometry.

\section{Preliminary Results in Sasakian geometry}
\setcounter{equation}{0}

\subsection{Transverse K\"ahler structure}

Let $(M, \xi, \eta , \Phi , g)$ be a $2m+1$-dimensional Sasakian
 manifold, and let $\mathcal F _{\xi }$ be the
characteristic foliation generated by $\xi $. Firstly, Let us recall that a Sasakian structure induce a transverse K\"ahler structure on the foliation $\mathcal F_{\xi }$. A transverse
holomorphic structure on $\mathcal F_{\xi }$ is given by an open
covering $\{U_{i}\}_{i\in A}$ of $M$ and local submersion $f_{i}:
U_{i}\rightarrow C^{n}$ with fibers of dimension $1$, (the leaves
of the foliation $\mathcal F _{\xi }|_{U_j}$ on $U_j$ coincide with the fibers of the map $f_{i}$, leave is the image of the flow of $\xi $), such that
for $i, j \in A$ there is a holomorphic isomorphism $\theta_{ij}$
of open sets of $C^{n}$ such that $f_{i}=\theta_{ij}\circ f_{j}$
on $U_{i}\cap U_{j}$.

\medskip

 In order to consider the deformations of
Sasakian structures, we consider the quotient bundle of the
foliation $\mathcal F _{\xi }$, $\nu (\mathcal F _{\xi })
=TM/L\xi$. The metric $g$ gives a bundle isomorphism $\sigma$
between $\nu (\mathcal F _{\xi })$ and the contact sub-bundle
$\mathcal D =Ker \eta$, where $\sigma : \nu (\mathcal F _{\xi })\rightarrow
\mathcal D$ defined by
\begin{eqnarray*}
\sigma ([X])=X-\eta (X)\xi .
\end{eqnarray*}
 By this isomorphism, $\Phi |_{\mathcal D}$ induces a complex
structure $\bar{J}$ on $\nu (\mathcal F _{\xi })$. Since the Nijenhuis torsion tensor
of $\Phi$ satisfies
\begin{eqnarray*}
N_{\Phi }(X , Y)=-d\eta (X , Y)\otimes \xi .
\end{eqnarray*}
So, $(\nu (\mathcal F_{\xi }) , \bar{J})\cong(\mathcal D,
\Phi |_{\mathcal D} )$ gives $\mathcal F _{\xi }$ a transverse
holomorphic structure. Then $(\mathcal D,
\Phi |_{\mathcal D} , d\eta)$ gives $\mathcal F_{\xi }$ a
 a transverse K\"ahler
structure with transverse K\"ahler form $\frac{1}{2}d\eta$ and
metric $g^{T}$ defined by $g^{T}(\cdot , \cdot )=\frac{1}{2}d\eta
(\cdot , \Phi \cdot)$.

In the following we say that a Sasakian structure $(\xi, \eta ' , \Phi ', g')$ have the same transverse holomorphic transverse structure with that of  $(\xi, \eta , \Phi , g)$, means that it satisfies
\begin{eqnarray}
 \bar{J}\circ \pi_{\nu } =\pi_{\nu } \circ \Phi '
\end{eqnarray}
where $\pi_{\nu}$ is the projection $\pi_{\nu}: TM\rightarrow \nu (\mathcal F_{\xi})$.

\medskip

{\bf Definition 2.1. } {\it We define $\mathcal S (\xi , \bar{J})$ to be the set of all Sasakian structures which have the same Reeb vector field and the same transverse holomorphic structure with $(\xi, \eta , \Phi , g)$, i.e. all Sasakian structures which compatible with $(\xi, \eta , \Phi , g)$.}

\medskip

{\bf Definition 2.2. } {\it Fixed a transverse holomorphic structure $(\nu (\mathcal F_{\xi }) , \bar{J})$ on  the characteristic foliation $\mathcal F_{\xi }$.  A complex vector field $X$ on  $M$ is called a  transverse holomorphic vector
field  if it satisfies:
\begin{enumerate}\item $\pi [\xi , X] =0$; \item $\bar{J}(\pi
(X))=\sqrt{-1}\pi (X)$; \item $\pi ([Y, X])-\sqrt{-1}\bar{J}\pi ([Y
, X])=0$, $\forall Y$ satisfying $\bar{J}\pi
(Y)=-\sqrt{-1}\pi (Y)$,\end{enumerate}
where $\pi $ is the projection $\pi : TM\rightarrow \nu (\mathcal F_{\xi })$.
Given a transverse K\"ahler form $d\eta $. Let $\psi $ be a complex valued basic function, then there is a unique vector field $V_{d\eta }(\psi )\in \Gamma (T^{c}M)$
satisfies: (1)$\bar{J}(\pi
(V_{d\eta }(\psi )))=\sqrt{-1}\pi (V_{d\eta }(\psi ))$; (2) $\psi =\sqrt{-1}\eta (V_{\eta }(\psi ))$; (3)
$\bar{\partial }_{B}\psi =-\frac{\sqrt{-1}}{2}d\eta (V_{\eta
}(\psi ) , \cdot )$. The vector field $V_{\eta }(\psi )$ is
called the  Hamiltonian vector field of $\psi $
corresponding to the transverse K\"ahler form $d\eta $. A complex vector field $X$ on  $M$ is called a  Hamiltonian holomorphic vector
field  if it is transverse holomorphic and is the Hamiltonian vector field of some complex valued basic function $\psi $
corresponding to some transverse K\"ahler form $d\eta $.}

\medskip

{\bf Remark 2.3. } {\it By the definition,  $c\xi $ is a Hamiltonian holomorphic vector field for any constant $c$. In this paper, without nontrivial Hamiltonian holomorphic vector
field  means that any Hamiltonian holomorphic vector
field must be $0$ or $c\xi $.}

\medskip

Given a Sasakian structure $(\xi , \eta , \Phi , g)$, we might identity $\mathcal D$ with $\nu (\mathcal F _{\xi })$ by the isomorphism. However, it's better to distinguish them, since under the deformations of Sasakian structure, the contact sub-bundle $\mathcal D$ changes, while $\nu (\mathcal F _{\xi })$ keeps fixed. But for simplicity of notation, we will use the same notation if there is no confusion, especially if we do not consider deformations. From the   transverse K\"ahler structure  $(\mathcal D,
\Phi |_{\mathcal D} , d\eta)$, one can
define the transverse Levi-Civita connection $\nabla^{T}$ on
$\mathcal D$ by
\begin{eqnarray}
\nabla^{T}_{X}Y=\left\{ \begin{array}{ll} &(\nabla_{X}Y)^{p},
\quad
 X\in \mathcal D,\\
 & [\xi , Y]^{p}, \quad  X=\xi , \\
\end{array}\right.
\end{eqnarray}
where $\nabla $ is the Levi-Civita connection with respect to the Riemannian metric $g$, $Y$ is a section of $\mathcal D$ and $X^{p}$ the projection
of $X$ onto $\mathcal D$.
It is easy to check that the
transverse Levi-Civita connection is torsion-free and metric
compatible. The transverse curvature tensor and transverse Ricci curvature are defined by
\begin{eqnarray}
R^{T}(V , W)Z =\nabla ^{T}_{V} \nabla ^{T}_{W} Z -\nabla ^{T}_{W}
\nabla ^{T}_{V} Z-\nabla ^{T}_{[V , W]}Z,
\end{eqnarray}
\begin{eqnarray}
Ric^{T}(X , Y) =<R^{T}(X , e_{i})e_{i} , Y>_{g},
\end{eqnarray}
where $e_{i}$ is an orthonormal basis of $\mathcal D$, $X , Y, Z \in \mathcal D$ and $V , W\in TM$.
We also have the following relations between the transverse
curvature tensor and the Riemann curvature tensor (see \cite{5})
\begin{eqnarray}\label{TRCR}
\begin{array}{lll}
&&R^{T}(X , Y)Z
=R(X , Y)Z -<R(X , Y)Z , \xi
>\xi \\
&&-<\nabla_{Y} Z , \xi
>\Phi (X)+<\nabla_{X}Z , \xi >\Phi (Y)
+<[X , Y] , \xi >\Phi (Z),\\
\end{array}
\end{eqnarray}
and
\begin{eqnarray}\label{2.14}
Ric^{T}(X, Y) =Ric (X, Y )+2g^{T}(X, Y),
\end{eqnarray}
for $X, Y, Z \in \mathcal D$.
The transverse Ricci form defined as following
\begin{eqnarray}
\rho ^{T}(X, Y)=Ric^{T}(\Phi X , Y).
\end{eqnarray}

\medskip

{\bf Definition 2.4. } {\it A Sasakian manifold $(M, \xi , \eta ,
\Phi , g)$ is said to be transversely K\"ahler-Einstein  if
\begin{eqnarray*}
Ric^{T}=\mu g^{T}, \quad or \quad \rho^{T}=\mu (\frac{1}{2}d\eta )
\end{eqnarray*}
for some constant $c$. }

\medskip

A Sasakian manifold $(M, \xi , \eta , \Phi , g)$ is called a
Sasakian-Einstein manifold if $g $  is a Riemannian Einstein
metric, i.e. $Ric_{g}=cg$ for some constant $c$. It is easy to see
that a
Sasakian-Einstein must be transversely K\"ahler-Einstein, the constant $c$ must be $2m$, and
\begin{eqnarray}\label{chern1}
Ric^{T}=(2m+2) g^{T}, \quad or \quad \rho^{T}=(m+1) d\eta  .
\end{eqnarray}

\medskip

\subsection{ Basic cohomology}

A $p$-form $\theta $ on Sasakian manifold $(M, \xi, \eta, \Phi,  g)$ is called
basic if
\begin{eqnarray}
i_{\xi }\theta =0, \quad L_{\xi }\theta =0,
\end{eqnarray}
where $i_{\xi}$ is the contraction with the Killing vector field
$\xi $, $L_{\xi}$ is the Lie derivative with respect to $\xi$.  Basic cohomology was introduced by Reinhart in \cite{[Rei59b]}. We begin with a brief review following \cite{[Ton97]}. It
is easy to see that the exterior differential preserves  basic
forms. Namely, if $\theta $ is a basic form, so is $d\theta $. Let
$\wedge _{B}^{p}(M)$ be the sheaf of germs of  basic $p$-forms and
$\Omega _{B}^{p}(M)=\Gamma (M , \wedge_{B}^{p}(M))$ the set of all
section of $\wedge _{B}^{p}(M)$. The  basic cohomology can be
defined in a usual way (see \cite{11}). Let $\mathcal D^{C}$ be
the complexification of the sub-bundle $\mathcal D$, and decompose
it into its eigenspaces with respect to $\Phi|_{D}$, that is
\begin{eqnarray}
\mathcal D^{C}= \mathcal D^{1 , 0} \oplus \mathcal D^{0 , 1}.
\end{eqnarray}
Similarly, we have a splitting of the complexification of the
bundle $\wedge ^{1}_{B}(M)$ of basic one forms on $M$,
\begin{eqnarray}
\wedge ^{1}_{B}(M) \otimes C =\wedge ^{1, 0}_{B}(M)\oplus \wedge
^{0, 1}_{B}(M).
\end{eqnarray}
Let $\wedge_{B} ^{i , j}(M) $ denote the bundle of basic forms of
type $(i, j)$.  Accordingly, we have the following decomposition
\begin{eqnarray}
\wedge _{B}^{p}(M)\otimes C =\oplus _{i+j=p}\wedge_{B} ^{i ,
j}(M).
\end{eqnarray}
Define $\partial_{B}$ and $\overline{\partial }_{B}$ by
\begin{eqnarray}
\begin{array}{lll}
\partial_{B} & : & \wedge_{B} ^{i ,
j}(M)\rightarrow \wedge_{B} ^{i+1 , j}(M);\\
\overline{\partial }_{B} & : & \wedge_{B} ^{i ,
j}(M)\rightarrow \wedge_{B} ^{i , j+1}(M);\\
\end{array}
\end{eqnarray}
which is the decomposition of $d$. Let
$d^{c}_{B}=\frac{1}{2}\sqrt{-1}(\overline{\partial
}_{B}-\partial_{B})$ and $d_{B}=d|_{\wedge^{p}_{B}}$. We have
$d_{B}=\overline{\partial }_{B}+\partial_{B}$,
$d_{B}d^{c}_{B}=\sqrt{-1}\partial_{B}\overline{\partial }_{B}$,
$d_{B}^{2}=(d^{c}_{B})^{2}=0$. The basic cohomology groups $H^{i,
j}_{B}(M, \mathcal F_{\xi})$ are fundamental invariants of a
Sasakian structure which enjoy many of the same properties as the
Dolbeault cohomology of a K\"ahler structure \cite{Boyer}. On
Sasakian manifolds, the $\partial \overline{\partial }$-lemma
holds for basic forms.

\medskip

{\bf Proposition 2.5. (\cite{11}) } {Let $\theta $ and $\theta '$
be two real closed basic form of type $(1, 1)$ on a compact
Sasakian manifold $(M, \xi , \eta , \Phi , g)$. If $[\theta
]_{B}=[\theta ']_{B}\in H^{1, 1}_{B}(M, \mathcal F_{\xi})$, then
there is a real basic function $\varphi $ such that
\begin{eqnarray*}
\theta =\theta ' +\sqrt{-1}\partial _{B} \overline{\partial
}_{B}\varphi .
\end{eqnarray*}

}

\medskip

Consider the complex bundle $(\mathcal D, \Phi |_{\mathcal D})$ (or $(\nu (\mathcal F_{\xi }), \bar{J})$) on a Sasakian manifold $(M, \xi , \eta, \Phi, g)$. Let $Rm^{T}$ be the transverse curvature with respect to  the transverse Levi-Civita connection $\nabla ^{T}$. If we choose a local foliate transverse frame $(X_{1}, \cdots , X_{m})$ on the bundle $\mathcal D$, then $Rm^{T}$ can be seen as a matrix valued 2-form (i.e. $End (\mathcal D)$-valued). $Rm^{T}$ is a basic $(1, 1)$-form. Let's define the basic $(k, k)$-form $\gamma_{k}$ by the formula
\begin{eqnarray}
det (Id_{m}+\frac{\sqrt{-1}}{2\pi }Rm^{T})=1+\sum_{k=1}^{m}\gamma_{k}.
\end{eqnarray}

\medskip

{\bf Definition 2.6. } {$\gamma_{k}$ is a closed basic $(k, k)$-form it represents an element in $H^{k, k}_{B}(M, \mathcal F_{\xi})$ that is called the basic $k^{th}$ Chern class and denoted by $c_{k}^{B}(M, \mathcal F_{\xi})$. }

\medskip

We have the following proposition.

\medskip

{\bf Proposition 2.7. (\cite{Boyer}, proposition 7.5.21)} {\it The basic Chern classes $c_{k}^{B}(M, \mathcal F_{\xi})$ are independent of the choice of a Sasakian structure in $\mathcal S (\xi, \bar{J})$}

\medskip

Let $\rho^{T}=Ric^{T}(\Phi \cdot , \cdot )$ be
the transverse Ricci form of the Sasakian structure $(\xi , \eta, \Phi, g)$.  $\rho^{T}$ is a real
closed basic $(1, 1)$-form and the basic cohomology class
$\frac{1}{2\pi}[\rho ^{T}]_{B}=c_{1}^{B}(M , \mathcal F _{\xi })$
is the basic first Chern class. We say that the basic first Chern
class of $(M, \mathcal F_{\xi })$ is positive (negative , null
resp.) if it contains a positive (negative, null resp. )
representation.
By the definition and (\ref{chern1}), a necessary condition for the existence of
Sasakian-Einstein metric $(M, \xi , \eta ,
\Phi , g)$ is $2\pi c_{1}^{B}(M, \mathcal F_{\xi })=(m+1)[d\eta]_{B}$.

\medskip

On Sasakian manifold $(M, \xi , \eta, \Phi, g)$, the basic Laplacian is defined by
\begin{eqnarray}
\triangle_{B} \psi =\frac{4m\sqrt{-1}\partial_{B}\bar{\partial}_{B}\psi \wedge (d\eta)^{m-1}\wedge \eta }{(d\eta)^{m}\wedge \eta },
\end{eqnarray}
for any basic function $\psi$. It is well known that the basic Laplacian is equal to the restriction of the Riemannian Laplacian $\triangle_{g} $ on
basic functions, i.e $\triangle_{B}\psi =\triangle_{g}\psi $ for any basic function $\psi$.

\subsection{Transformations of Sasakian structures} Let $(\xi , \eta , \Phi , g)$ be a Sasakian structure on $M$, for every real basic function $\varphi \in \mathcal H (\xi , \eta , \Phi , g)$, we can obtain a new Sasakian structure $(\xi , \eta_{\varphi } , \Phi_{\varphi} , g_{\varphi })$.  Here $\eta_{\varphi } , \Phi_{\varphi} , g_{\varphi }$ is defined as that in (\ref{eta-d}) and (\ref{phi-d}). We know that the above deformations  fix the Reeb vector field $\xi $ and the transverse holomorphic
structure.

Let $(\xi ' , \eta ' , \Phi ', g ' ) $ be an another Sasakian structure
which is compatible with
$(\xi  , \eta  , \Phi  , g )$,
 then  there exists a basic function $\varphi \in \mathcal H (\xi , \eta , \Phi , g)$ such that
$
d\eta '=d\eta +dd^{c}_{B}\varphi$ and $\eta '=\eta +d_{B}^{c}\varphi +\zeta
$,
where $\zeta $ is closed basic one form. Since the difference between $\eta '$ and $\eta $ is a basic one form, it is easy to see that
\begin{eqnarray}
(d\eta ')^{m}\wedge \eta '=(d\eta ')^{m}\wedge \eta ,
\end{eqnarray}
and
\begin{eqnarray}
\int_{M}(d\eta ')^{m}\wedge \eta '=\int_{M}(d\eta )^{m}\wedge \eta =V.
\end{eqnarray}
In the following, $\rho ^{T}_{d\eta }$ denotes the transverse Ricci form with respect to the transverse K\"ahler metric $d\eta$. In local foliation coordinates $(x, z_{1}, \cdots , z_{m})$, we have
\begin{eqnarray}
\rho ^{T}_{d\eta }=-\sqrt{-1}\partial_{B}\overline{\partial
}_{B} \log det ((d\eta)_{i\bar{j}}).
\end{eqnarray}
Globally, the difference between two transverse Ricci form can be expressed by
\begin{eqnarray}
\rho ^{T}_{d\eta }-\rho ^{T}_{d\eta '}=\sqrt{-1}\partial_{B}\overline{\partial }_{B}\log
(\frac{(d\eta ')^{m}\wedge \eta '}{(d\eta )^{m}\wedge
\eta}).
\end{eqnarray}
From the above formula, we know that to find a Sasakian-Einstein structure compatible with $(\xi ,\eta, \Phi, g)$ is equivalent to solve the transverse Monge-Amp\`ere equation (\ref{MA1}).

\medskip

Let's recall one class special deformations of Sasakian structure,
\begin{eqnarray}
\xi_{s}=s^{-1}\xi, \quad \eta_{s}=s\eta , \quad \Phi_{s}=\Phi ,\quad
g_{s}=sg+s(s-1)\eta \otimes \eta,
\end{eqnarray}
where $s$ is constant. These were called  $\mathcal D$-homothetic
deformations by Tanno \cite{Tan}. These deformations do not deform
the characteristic foliation and the contact bundle $\mathcal D$,
but only rescale the Reeb field $\xi$ and contact $1$-form $\eta$.
If $(\xi , \eta , \Phi , g)$ be a transverse K\"ahler Einstein
Sasakian metric with positive transverse Ricci curvature,  i.e.
$Ric^{T}_{g}=\mu g^{T}$ ($Ric =(\mu -2 )g +(2n+2 -\mu ) \eta \otimes \eta $)for positive constant $\mu \in R$. Then,
by the $\mathcal D$-homothetic transformation $s =\frac{\mu}{2(n+1)}$, we
get a Sasakian-Einstein metric $(\xi_{s}, \eta_{s} , \Phi_{s} ,
g_{s}) $.
Indeed, by the relation formula (2.15) in \cite{Tan} , we have
\begin{eqnarray}
\begin{array}{lll}
Ric_{g_{s}}&=& Ric_{g}-2(s-1)g+s^{-1}(s-1)(2(2n+1)s +2ns(s-1))\eta \otimes \eta \\
&=&\{2n+2-\mu +s^{-1}(s-1)(2(2n+1)s +2ns(s-1))\} \eta \otimes \eta\\
&& +(\mu -2 -2(s-1))g \\
&=& \frac{\mu -2s}{s} \{sg +s(s-1)\eta \otimes \eta \} +((2n+2)s^{2}-\mu s )\eta \otimes \eta \\
&=& 2m g_{s}.\\
\end{array}
\end{eqnarray}

Through Tanno's $\mathcal D$-homothetic transformation, one can prove the following lemma. The proof can be found in \cite{Ni} and
\cite{Se}, see also proposition 2.6 in \cite{Z2}.

\medskip

{\bf Proposition 2.7. } {\it Let $(M, \xi , \eta , g)$ be a
$2m+1$-dimensional compact Sasakian manifold with $Ric^{T}\geq \epsilon g^{T}$.
Suppose that $\phi \in C_{B}^{\infty}(M)$ satisfies
$\triangle_{g}\phi \leq \delta $, then we have
\begin{eqnarray}
-\inf_{M} \phi \leq \frac{1}{V}\int_{M}(-\phi )(d\eta )^{m}\wedge
\eta + \frac{C(m)\delta }{\epsilon},
\end{eqnarray}
where $V=\int_{M} (d\eta )^{m}\wedge \eta $ and constant $C(m)$
depends only on $m$.
 }

\medskip

From the above proposition, it is easy to conclude the following corollary.

\medskip

{\bf Corollary 2.8. } {\it Let $(M, \xi , \eta , g)$ be a
$2m+1$-dimensional compact Sasakian manifold and $\varphi \in \mathcal H (\xi , \eta , \Phi , g)$ be a potential function with $Ric^{T}(d\eta_{\varphi})\geq \epsilon g^{T}_{\varphi }$.
Then we have
\begin{eqnarray}
Osc (\varphi ) \leq I_{d\eta }(\varphi ) + \frac{\tilde{C}(m) }{\epsilon} +\tilde{C}(M, g),
\end{eqnarray}
where constant $\tilde{C}(M, g) $ depends only on $(M, g)$ and constant $\tilde{C}(m)$
depends only on $m$. To the definition of $I_{d\eta }$, see section 3.
 }

\medskip

{\bf Proof. } Using the fact $\triangle_{d\eta }\varphi_{t}\geq -4m$ and the Green's formula, we have
\begin{eqnarray}
\sup_{M}\varphi\leq \frac{1}{V}\int_{M}\varphi(d\eta
)^{m}\wedge \eta +\tilde{C}(M, g).
\end{eqnarray}
Then, by
$\triangle_{d\eta_{\varphi }}\varphi\leq 4m$ and the above proposition, we have
\begin{eqnarray}
\begin{array}{lll}
Osc (\varphi )&= & \sup_{M}\varphi_{t}- \inf_{M}\varphi_{t} \\
&\leq & I_{d\eta }(\varphi ) + \frac{\tilde{C}(m) }{\epsilon} +\tilde{C}(M, g).
\end{array}
\end{eqnarray}

\hfill $\Box$

\section{Energy functionals}
\setcounter{equation}{0}

Let $( \xi , \eta , \Phi , g)$ be a Sasakian structure on $M$. We consider following  functionals on
$\mathcal H( \xi , \eta , \Phi , g)$ which are analogous to the ones
 in K\"ahler geometry.
\begin{eqnarray}\label{DF}
\begin{array}{lll}
I_{d\eta }(\varphi ) &:=& \frac{1}{V}\int_{M}\varphi \{(d\eta
)^{m}\wedge \eta -(d\eta_{\varphi } )^{m}\wedge \eta
\}\\
J_{d\eta } (\varphi ) &:= & \int_{0}^{1}\frac{1}{s}I_{d\eta
}(s\varphi ) ds,\\
F_{d\eta }^{0} (\varphi ) &:= & J_{d\eta } (\varphi )- \frac{1}{V}\int_{M}\varphi (d\eta
)^{m}\wedge \eta \\
F_{d\eta } (\varphi ) &:= & F_{d\eta }^{0} (\varphi )-\frac{1}{m+1}\log \{\frac{1}{V}\int_{M}e^{h_{d\eta }-(m+1)\varphi }(d\eta
)^{m}\wedge \eta\}\\
\end{array}
\end{eqnarray}
where $V=\int_{M}(d\eta )^{m}\wedge \eta $,  and $h_{d\eta}$ is a smooth basic function
which satisfies $\rho^{T}=(m+1)d\eta
+\sqrt{-1}\partial_{B}\bar{\partial }_{B}h_{d\eta}$ and $\int_{M}\exp(h_{d\eta}) (d\eta )^{m}\wedge \eta =V$. Noting that $g(\Phi
\cdot , \cdot )=\frac{1}{2}d\eta $, when $\varphi \in
C_{B}^{\infty}(M)$ we have
\begin{eqnarray}
m\sqrt{-1} \partial_{B}\bar{\partial }_{B}\varphi \wedge (d\eta
)^{m-1}\wedge \eta =\frac{1}{4}\triangle \varphi (d\eta
)^{m}\wedge \eta ,
\end{eqnarray}
where $\triangle $ is the Laplace of the metric $g$. Let
$\varphi_{s}$ be a smooth curve in $\mathcal H$, by direct
calculation, we have
\begin{eqnarray}
\begin{array}{lll}
\frac{d}{d s}I_{d\eta
}(\varphi_{s})&=&\frac{1}{V}\int_{M}\dot{\varphi}_{s}\{(d\eta
)^{m}-(d\eta_{\varphi_{s}})^{m}\}\wedge \eta \\
&&-\frac{1}{4V}\int_{M}\varphi_{s}\triangle_{\varphi_{s}}\dot{\varphi}_{s}
(d\eta_{\varphi_{s}})^{m}\wedge \eta ,
\end{array}
\end{eqnarray}
\begin{eqnarray}\label{DJ}
\frac{d}{d s}J_{d\eta
}(\varphi_{s})=\frac{1}{V}\int_{M}\dot{\varphi}_{s}\{(d\eta
)^{m}-(d\eta_{\varphi_{s}})^{m}\}\wedge \eta ,
\end{eqnarray}
\begin{eqnarray}\label{F0}
\frac{d}{ds}
F_{d\eta }^{0} (\varphi_{s} ) = -\frac{1}{V}\int_{M}\dot{\varphi}_{s}(d\eta_{\varphi_{s}})^{m}\wedge \eta ,
\end{eqnarray}
and
\begin{eqnarray}\label{F}
\begin{array}{lll}
&&\frac{d}{ds}
F_{d\eta } (\varphi_{s} ) = -\frac{1}{V}\int_{M}\dot{\varphi}_{s}(d\eta_{\varphi_{s}})^{m}\wedge \eta \\
 &&+(\int_{M}e^{h_{d\eta }-(m+1)\varphi }(d\eta
)^{m}\wedge \eta)^{-1}\int_{M}\dot{\varphi}_{s}e^{h_{d\eta }-(m+1)\varphi }(d\eta
)^{m}\wedge \eta , \\
\end{array}
\end{eqnarray}
where $\dot{\varphi}_{s}=\frac{d}{d s}\varphi_{s}$ and
$\triangle_{\varphi_{s}}$ is the Laplace corresponding with the
metric $g_{\varphi_{s}}$. From (\ref{F}), it is easy to that the critical points of $F_{d\eta}$ are transverse K\"ahler-Einstein metrics.

The following properties can be proved by a similar method as
in the K\"ahler case (see \cite{BM}, \cite{T1}, \cite{Ni},
\cite{Se}).

\medskip

{\bf Proposition 3.1. } {\it Let $C$ be a constant, then
\begin{eqnarray}
I_{d\eta }(\varphi +C)=I_{d\eta }(\varphi ), \quad J_{d\eta }(\varphi
+C)=J_{d\eta }(\varphi ), \quad F_{d\eta }(\varphi +C)=F_{d\eta }(\varphi ).
\end{eqnarray}
$I_{d\eta }$, $I_{d\eta }-J_{d\eta }$, $J_{d\eta }$ are non-negative
functionals on $\mathcal H( \xi , \eta , \Phi , g)$, and we have
\begin{eqnarray}\label{2.17}
I_{d\eta }(\varphi )\leq (m+1)\{I_{d\eta }(\varphi )-J_{d\eta
}(\varphi )\}\leq mI_{d\eta }(\varphi ).
\end{eqnarray}
Let ${\varphi_{t}}$ be a family of basic functions in $\mathcal
H$, then
\begin{eqnarray}\label{2D}
\frac{d}{dt}\{I_{d\eta }(\varphi_{t} )-J_{d\eta }(\varphi_{t}
)\}=\frac{-1}{4V}\int_{M}\varphi_{t}(\triangle_{t}
\frac{d}{dt}\varphi_{t})(d\eta_{\varphi_{t}})^{m}\wedge \eta ,
\end{eqnarray}
where $\triangle_{t}$ is the Laplacian  corresponding to the
metric $g_{\varphi_{t}}$. $F_{d\eta}$ satisfies the following cocycle property, i.e.
\begin{eqnarray}
F_{d\eta }(\psi )+F_{d\eta ' }(\phi -\psi )=F_{d\eta }(\phi ),
\end{eqnarray}
and
\begin{eqnarray}
F_{d\eta }(\psi )=-F_{d\eta '}(-\psi )
\end{eqnarray}
for all $\phi , \psi \in \mathcal H( \xi , \eta , \Phi , g)$ and $d\eta '=d\eta +\sqrt{-1}\partial_{B}\overline{\partial }_{B}\psi$. We also have the cocycle condition for $F_{d\eta}^{0}$.}

\medskip

{\bf Lemma 3.2. } {\it Let $(\xi , \eta , \Phi , g)$ and $(\xi , \eta ' , \Phi ' , g')$ are two Sasakian structures with the same transverse holomorphic structure on $M$, and assume that $d\eta =d\eta ' +\sqrt{-1}\partial_{B}\bar{\partial }_{B}\phi $ for some basic function $\phi$. Then, we have
\begin{eqnarray}\label{UU1}
|I_{d\eta ' }(\varphi +\phi ) -I_{d\eta }(\varphi )|\leq (m+1) Osc (\phi )
\end{eqnarray}
for all $\varphi \in \mathcal H (\xi , \eta , \Phi , g)$.

}

{\bf Proof. }
By definition, we have
\begin{eqnarray}\label{UU2}
\begin{array}{lll}
&&V(I_{d\eta ' }(\varphi +\phi ) -I_{d\eta }(\varphi ))\\&=&\int_{M}\varphi ((d\eta ')^{m}-(d\eta )^{m})\wedge \eta +\int_{M}\phi ((d\eta ')^{m}-(d\eta _{\varphi })^{m})\wedge \eta ,\\
\end{array}
\end{eqnarray}
where $d\eta_{\varphi }=d\eta  +\sqrt{-1}\partial_{B}\bar{\partial }_{B}\varphi $. On the other hand, By direct calculation, we have
\begin{eqnarray}\label{UU3}
\begin{array}{lll}
&&|\int_{M}\varphi ((d\eta ')^{m}-(d\eta )^{m})\wedge \eta |\\
&=&|\int_{M}\varphi (-\sqrt{-1}\partial_{B}\bar{\partial }_{B}\phi)\wedge (\sum _{j=0}^{m-1}(d\eta ')^{j}\wedge (d\eta )^{m-j-1})\wedge \eta |\\
&=&|\int_{M}\phi (-\sqrt{-1}\partial_{B}\bar{\partial }_{B}\varphi)\wedge (\sum _{j=0}^{m-1}(d\eta ')^{j}\wedge (d\eta )^{m-j-1})\wedge \eta |\\
&=&|\int_{M}\phi (d\eta -d\eta_{\varphi } )\wedge (\sum _{j=0}^{m-1}(d\eta ')^{j}\wedge (d\eta )^{m-j-1})\wedge \eta | \\
&\leq & m V Osc (\phi )
\end{array}
\end{eqnarray}
and
\begin{eqnarray}\label{UU4}
|\int_{M}\phi ((d\eta ')^{m}-(d\eta _{\varphi })^{m})\wedge \eta|\leq V Osc (\phi ).
\end{eqnarray}
Then (\ref{UU2}), (\ref{UU3}) and (\ref{UU4}) imply (\ref{UU1}).

\hfill $\Box$

\medskip

Let $\rho^{T}_{\varphi }$ denote the transverse Ricci form of
the Sasakian structure  $(\xi , \eta _{\varphi } , \Phi_{\varphi
}, g_{\varphi }) $. $\int_{M} \rho^{T}_{\varphi } \wedge (d\eta_{\varphi })^{m}\wedge
\eta_{\varphi }$ is independent of the choice of $\varphi \in \mathcal
H (\xi , \eta , \Phi , g)$ (e.g., Proposition 4.4 in \cite{12}).
This means that
\begin{eqnarray}\label{scalar1}
\bar{S}=\frac{\int_{M} S_{\varphi}^{T}  (d\eta
_{\varphi})^{m}\wedge \eta_{\varphi}}{\int_{M} (d\eta_{\varphi}
)^{m}\wedge \eta_{\varphi} }=\frac{\int_{M} 2m \rho_{\varphi}^{T}
\wedge (d\eta_{\varphi} )^{m-1}\wedge \eta}{\int_{M}
(d\eta_{\varphi} )^{m}\wedge \eta },
\end{eqnarray}
depends only on the basic K\"ahler class.
As in the K\"ahler case (see \cite{mabuchi}), we can define the Mabuchi's $\mathcal K$-energy on the space $\mathcal H (\xi, \eta, \Phi, g)$.

\medskip

{\bf Definition 3.3. } {\it Let $\varphi '$ and $\varphi ''$ are two
basic functions in $\mathcal H (\xi, \eta, \Phi, g)$, we define
\begin{eqnarray}\label{Mabu}
\mathcal M (\varphi ' , \varphi ''):=-\frac{1}{V}\int_{0}^{1} \int_{M}
\dot{\varphi }_{t} (S ^{T}_{t}-\bar{S})(d\eta _{t})^{m}\wedge
\eta \, dt ,
\end{eqnarray}
where  $\varphi _{t}$ ($t\in [0, 1]$)
be a path in $\mathcal H$ connecting $\varphi '$ and $\varphi ''$, $\dot{\varphi
}_{t}=\frac{\partial }{\partial t}\varphi_{t}$, $S^{T}_{t}$ is the
transverse scalar curvature to the Sasakian structure $(\xi ,
\eta_{\varphi_t} , \Phi_{\varphi_t} , g_{\varphi_t})$ and $\bar{S}$ is the average defined
as in (\ref{scalar1}). We also define
\begin{eqnarray}
\mathcal V_{d\eta_{\varphi '}} (\varphi  ):=\mathcal M(\varphi ' , \varphi ' + \varphi )
\end{eqnarray}
for any $\varphi \in \mathcal H (\xi , \eta_{\varphi '} , \Phi _{\varphi '}, g_{\varphi '})$.

}

\medskip

By Theorem 4.12 in \cite{12} (or lemma 11 in \cite{GZ}), we know that $\mathcal M$
is independent of the path $\varphi _{t}$, and so $\mathcal M$ is well defined. Furthermore, $\mathcal M$ satisfies the
following cocycle condition, i.e.
\begin{eqnarray}
\mathcal M(\varphi _{0}, \varphi_{1})+\mathcal M(\varphi _{1},
\varphi_{0})=0,
\end{eqnarray}
\begin{eqnarray}
\mathcal M(\varphi _{0}, \varphi_{1})+\mathcal M(\varphi _{1},
\varphi_{2})+\mathcal M(\varphi _{2}, \varphi_{0})=0.
\end{eqnarray}
 and
\begin{eqnarray}\label{2.32}
\mathcal M (\varphi_{1}  +C' , \varphi _{2}+C'')=\mathcal M (\varphi _{1}
, \varphi _{2})
\end{eqnarray}
for any $\varphi_{i}\in \mathcal H (\xi, \eta, \Phi, g)$ and $C' , C'' \in \mathbb R$. Following Ding \cite{D1}, we also have the following relation between the functionals $\mathcal V_{d\eta }$ and $F_{d\eta}$.

\medskip

{\bf Remark: } {\it By the definitions and the above properties, it is easy to see that the functionals $F_{d\eta }$ and $\mathcal V_{d\eta }$ can also be defined on the space $\mathcal S(\xi , \bar{J})$.}

\medskip

{\bf Lemma 3.4. } {Let  $(M, \xi , \eta , \Phi , g)$ be a compact Sasakian manifold with $2\pi
c_{1}^{B}(M)=(m+1)[d\eta ]_{B}$, then
\begin{eqnarray}\label{MF1}
\begin{array}{lll}
&&\mathcal V_{d\eta} (\phi  )-2(m+1) F_{d\eta } (\phi )\\&=&\frac{2}{V}\int_{M} h_{d\eta }(d\eta )^{m}\wedge \eta -\frac{2}{V}\int_{M} h_{d\eta_{\phi }}(d\eta_{\phi} )^{m}\wedge \eta \\
\end{array}
\end{eqnarray}
for any $\phi \in \mathcal H (\xi, \eta, \Phi, g)$, where $h_{d\eta }$
and $h_{d\eta_{\phi }}$ are the normalized  Ricci potential functions with respect to $d\eta$ and $d\eta_{\phi}$.

}

\medskip

{\bf Proof. } Let $\phi_{t}$ be a path connecting $0$ with $\phi$, by the definition, we have
\begin{eqnarray}\label{MF2}
\begin{array}{lll}
\mathcal V_{d\eta} (\phi  )&=&-\frac{1}{V}\int_{0}^{1} \int_{M}
\dot{\phi }_{t} (S ^{T}_{t}-2m(m+1))(d\eta _{t})^{m}\wedge
\eta \, dt \\
&=&-\frac{2m}{V}\int_{0}^{1} \int_{M}
\dot{\phi }_{t} (\rho ^{T}_{t}-(m+1)d\eta_{t})(d\eta _{t})^{m-1}\wedge
\eta \, dt\\
&=&-\frac{2m\sqrt{-1}}{V}\int_{0}^{1} \int_{M}
\dot{\phi }_{t} \partial_{B}\bar{\partial }_{B}(h_{d\eta } -(m+1)\phi_{t}\\
&&-\log \frac{(d\eta _{t})^{m}\wedge
\eta}{(d\eta )^{m}\wedge
\eta})(d\eta _{t})^{m-1}\wedge
\eta \, dt\\
&=&-\frac{2}{V}\int_{0}^{1} \int_{M}
(h_{d\eta } -(m+1)\phi_{t}-\log \frac{(d\eta _{t})^{m}\wedge
\eta}{(d\eta )^{m}\wedge
\eta})\frac{\partial }{\partial t}((d\eta _{t})^{m}\wedge
\eta ) \, dt\\
&=& -2(m+1)(I_{d\eta }-J_{d\eta})(\phi ) +\frac{2}{V}\int_{M} h_{d\eta }(d\eta ^{m}-d\eta_{\phi}^{m})\wedge \eta  \\
&&+\frac{2}{V}\int_{M} \log \frac{(d\eta _{\phi})^{m}\wedge
\eta}{(d\eta )^{m}\wedge
\eta}(d\eta_{\phi})^{m}\wedge \eta .
\end{array}
\end{eqnarray}
On the other hand, it is easy to check that
\begin{eqnarray}
-\log \frac{(d\eta _{\phi})^{m}\wedge
\eta}{(d\eta )^{m}\wedge
\eta}-(m+1)\phi +c =h_{d\eta_{\phi }}-h_{d\eta},
\end{eqnarray}
where $c=-\log (\frac{1}{V}\int_{M}e^{h_{d\eta }-(m+1)\phi } d\eta^{m}\wedge \eta)$. Then
\begin{eqnarray*}
\begin{array}{lll}
&&\mathcal V_{d\eta} (\phi  )=2(m+1) J_{d\eta } (\phi )-\frac{2(m+1)}{V}\int_{M}\phi d\eta^{m}\wedge \eta +2c\\&&+\frac{2}{V}\int_{M} h_{d\eta }(d\eta )^{m}\wedge \eta -\frac{2}{V}\int_{M} h_{d\eta_{\phi }}(d\eta_{\phi} )^{m}\wedge \eta \\
&=&2(m+1) F_{d\eta } (\phi )+\frac{2}{V}\int_{M} h_{d\eta }(d\eta )^{m}\wedge \eta -\frac{2}{V}\int_{M} h_{d\eta_{\phi }}(d\eta_{\phi} )^{m}\wedge \eta . \\
\end{array}
\end{eqnarray*}

\hfill $\Box$

By the normalization $\int_{M}\exp(h_{d\eta_{\phi}}) (d\eta_{\phi} )^{m}\wedge \eta =V$, we known that $\int_{M}h_{d\eta_{\phi}} (d\eta_{\phi} )^{m}\wedge \eta \leq 0$, then we have the following corollary.

\medskip

{\bf Corollary 3.5. } {Let  $(M, \xi , \eta , \Phi , g)$ be a compact Sasakian manifold with $2\pi
c_{1}^{B}(M)=(m+1)[d\eta ]_{B}$, then
\begin{eqnarray}\label{MF3}
\mathcal V_{d\eta} (\phi  )\geq 2(m+1) F_{d\eta } (\phi )+\frac{2}{V}\int_{M} h_{d\eta }(d\eta )^{m}\wedge \eta
\end{eqnarray}
for any $\phi \in \mathcal H (\xi, \eta, \Phi, g)$.

}

 \hspace{0.4cm}

\section{transverse Monge-Amp\`ere equation}
\setcounter{equation}{0}

Let $(M, \xi , \eta , \Phi , g)$ be a compact Sasakian manifold with $2\pi
c_{1}^{B}(M)=(m+1)[d\eta ]_{B}$.  Given any $\varphi \in \mathcal H (M, \xi , \eta , \Phi ,
g)$, we have a new Sasakian structure $(\xi , \eta_{\varphi },
\Phi_{\varphi }, g_{\varphi })$.
It is easy to check that $d\eta_{\varphi }$ is Sasakian-Einstein if and only if
$\varphi $ satisfies the following equation
\begin{eqnarray}
\sqrt{-1}\partial_{B}\overline{\partial }_{B}\log
(\frac{(d\eta_{\varphi })^{m}\wedge \eta}{(d\eta )^{m}\wedge
\eta})=\sqrt{-1}\partial_{B}\overline{\partial
}_{B}(h_{d\eta }-(m+1)\varphi ),
\end{eqnarray}
which is equivalent to the  transverse Monge-Amp\`ere equation
(\ref{MA1}).  As in K\"ahler case, we consider a family of
equations (\ref{MA2}). We set
\begin{eqnarray}
S=\{t\in [0, 1]| (\ref{MA2})\quad is \quad solvable \quad for
\quad t\}.
\end{eqnarray}
By \cite{11}, we know that (\ref{MA2}) is solvable for $t=0$. The openness of $S$ was proved  in \cite{Se}(Proposition 5.3) and \cite{Ni}(Proposition
4.4) in a similar way as that in \cite{Au}.
$S$ is not empty. In order to use the continuity method to solve
(\ref{MA1}), we only need to prove the closedness of $S$. By El-Kacimi's ( \cite{11})
generalization of Yau's estimate (\cite{Y}) for transverse
Monge-Amp\`ere equations, the $C^{0}$-estimate for solutions of
(\ref{MA2}) implies the $C^{2, \alpha }$-estimate for them, and
the transverse elliptic Schauder estimates give higher order
estimates. Therefore it suffices to estimate $C^{0}$-norms of the
solutions of (\ref{MA2}). We list the following proposition for further discussion, the proof can be found in \cite{Se}(proposition 5.3) and \cite{Ni}(proposition 4.4), see also proposition 4.2 in \cite{Z2}.

\medskip

{\bf Proposition 4.1. } {\it Let $0< \tau \leq 1$, and suppose that
(\ref{MA2}) has a solution $\varphi _{\tau}$ at $t=\tau$. If $0<\tau <1$, then there exists
some $\epsilon >0$ such that $\varphi_{\tau }$ uniquely extends to a smooth
family of solution $\{\varphi _{t}\}$ of (\ref{MA2}) for $t\in (0,
1)\cap (\tau -\epsilon , \tau +\epsilon )$. $S$ is also open near
$t=0$, i.e. there exists a small positive number $\epsilon $ such that there is  a smooth  family solution of (\ref{MA2}) for $t\in (0, \epsilon)$. If $(M, \xi, \eta, \Phi, g)$ admits no nontrivial Hamiltonian holomorphic vector field, $\varphi_{1}$ can also be extended uniquely to a smooth
family of solution $\{\varphi _{t}\}$ of (\ref{MA2}) for $t\in (1-\epsilon , 1]$.}

\medskip

As in \cite{BM}, we have the following lemma, the proof can be found in \cite{Z2}
(lemma 4.3), see also lemma 4.9 in \cite{Ni} and lemma 5.4 in \cite{Se}.

\medskip

{\bf Lemma 4.2. } {Let  $\{\varphi_{t}\}$ be a smooth family
of solution  of (\ref{MA2}) for $t\in (0, 1]$, then
\begin{eqnarray}\label{mono1}
\frac{d}{dt}(I_{d\eta }-J_{d\eta })(\varphi_{t})\geq 0.
\end{eqnarray}

}

\medskip

Now, we consider the existence problem of Sasakian-Einstein
metrics. We prove the following theorem.

\medskip

{\bf Theorem 4.3. } {\it Let  $(M, \xi , \eta , \Phi , g)$ be a compact Sasakian manifold with $[d\eta ]_{B} = \frac{2\pi }{m+1}
c_{1}^{B}(M, \mathcal F_{\xi})$. If $F_{d\eta }$ (or $\mathcal V_{d\eta}$) is proper in the space $\mathcal H(\xi , \eta , \Phi , g)$, then there must exists a Sasakian-Einstein metric compatible with $(\xi , \eta , \Phi , g)$ .}

\medskip

{\bf Proof. }
By
Proposition 4.1, we can suppose that there exists a smooth family
of solution $\{\varphi_{t}\}$ of (\ref{MA2}) for $t\in (0, \tau )$
with some $\tau \in (0, 1)$.
 From the equation (\ref{MA2}), we know that $\triangle_{t}\varphi_{t}\leq 4m$ and $\rho^{T}_{d\eta_{t}}\geq t(m+1)d\eta_{t}$. By proposition 2.7, we have
\begin{eqnarray}\label{G2}
\frac{1}{V}\int_{M}\varphi_{t}(d\eta_{\varphi}
)^{m}\wedge \eta \leq inf_{M}\varphi_{t}+\frac{C_{1}(m)}{t}.
\end{eqnarray}
where positive constant $C_{1}(m)$ depends only on $m$.
 Using the fact $\triangle_{d\eta }\varphi_{t}\geq -4m$ and the Green formula, we have
\begin{eqnarray}
\sup_{M}\varphi_{t}\leq \frac{1}{V}\int_{M}\varphi_{t}(d\eta
)^{m}\wedge \eta +C_{2}
\end{eqnarray}
where $C_{2}$ is a positive constant depends only on the geometry of $(M, g)$. By the normalization, it's easy to check that  $\sup_{M}\varphi_{t} \geq 0$ and $\inf_{M}\varphi_{t}\leq 0$. Then
\begin{eqnarray}\label{C0}
\begin{array}{lll}
\|\varphi_{t}\|_{C^{0}}&\leq & \sup_{M}\varphi_{t}- \inf_{M}\varphi_{t} \\
&\leq & I_{d\eta }(\varphi_{t}) +\frac{C_{1}(m)}{t}+C_{2}.
\end{array}
\end{eqnarray}
By (\ref{2.17}) and (\ref{mono1}), we have
\begin{eqnarray}\label{mono2}
I_{d\eta}(\varphi_{t_{1}})\leq (m+1)(I_{d\eta }-J_{d\eta })(\varphi_{t_{2}})
\end{eqnarray}
for any $0<t_{1}\leq t_{2}<\tau$. Combining (\ref{C0}) and (\ref{mono2}), we get
\begin{eqnarray}
t\|\varphi_{t}\|_{C^{0}}\leq t_{0}(m+1)(I_{d\eta }-J_{d\eta })(\varphi_{t_{0}})+C_{3}
\end{eqnarray}
for any $0<t\leq t_{0}<\tau$, where $C_{3}$ is a positive constant depends only on the geometry of $(M, g)$. So, we obtain an uniform bound on
$|\frac{(d\eta+\sqrt{-1}\partial_{B}\overline{\partial }_{B}\varphi_{t}
)^{m}\wedge \eta}{(d\eta )^{m}\wedge \eta}|$ for $0<t\leq t_{0}<\tau$. By El-Kacimi's ( \cite{11})
generalization of Yau's $C^{0}$ estimate (\cite{Y}) for transverse
Monge-Amp\`ere equations, there exists an uniform constant $C_{4}$ such that
\begin{eqnarray}\label{C02}
\|\varphi_{t}\|_{C^{0}}\leq C_{4}
\end{eqnarray}
for $0<t\leq t_{0}<\tau$.

Differentiating (\ref{MA2}) with respect to $t$, we have
\begin{eqnarray}\label{DM}
\frac{1}{4}\triangle_{t}\dot{\varphi}_{t}=-t(m+1)\dot{\varphi}_{s} -(m+1)\varphi_{t}.
\end{eqnarray}
Using (\ref{F0}) and (\ref{DM}), we have
\begin{eqnarray}\label{FFF1}
\begin{array}{lll}
\frac{d}{dt} (tF_{d\eta }^{0}(\varphi_{t})) &=& F_{d\eta }^{0}(\varphi_{t}) -\frac{t}{V}\int_{M}\dot{\varphi}_{t} (d\eta_{\varphi })^{m}\wedge \eta\\
&=& J_{d\eta }(\varphi_{t})-\frac{1}{V}\int_{M}\varphi_{t} (d\eta)^{m}\wedge \eta -\frac{t}{V}\int_{M}\dot{\varphi}_{t} (d\eta_{\varphi })^{m}\wedge \eta\\
&=&-( I_{d\eta }(\varphi_{t})- J_{d\eta }(\varphi_{t}))\leq 0.\\
\end{array}
\end{eqnarray}
By the uniform $C^{0}$ estimate (\ref{C02}), it is easy to check that
\begin{eqnarray}
tF_{d\eta }^{0}(\varphi_{t}) \rightarrow 0
\end{eqnarray}
as $t\rightarrow 0$. So, from (\ref{FFF1}), we have
\begin{eqnarray}\label{F1}
\begin{array}{lll}
F_{d\eta } (\varphi_{t} )&\leq & -\frac{1}{m+1}\log \{\frac{1}{V}\int_{M}e^{h_{d\eta }-(m+1)\varphi }(d\eta
)^{m}\wedge \eta\}\\
&\leq & \frac{(1-t)}{V}\int_{M}\varphi_{t}(d\eta_{\varphi}
)^{m}\wedge \eta,\\
\end{array}
\end{eqnarray}
where we have used the concavity of the logrithmic function.
From (\ref{F1}) and (\ref{G2}), we have
\begin{eqnarray}\label{F2}
F_{d\eta } (\varphi_{t} )\leq \frac{(1-t)}{t}C_{1}(m).
\end{eqnarray}
By (\ref{F2}) and (\ref{C0}), the properness of $F_{d\eta }$ implies that $J_{d\eta }(\varphi_{t})$, and consequently, $\|\varphi_{t}\|_{C^{0}}$ is uniformly bounded for $t\in [\epsilon , \tau )$. Therefore, the equation (\ref{MA1}) can be solved, i.e. there is a Sasakian-Einstein metric on $M$.

\medskip

For the $\mathcal K$-energy case. It is easy to see that, along the solutions of (\ref{MA2}), we have
\begin{eqnarray}
S_{t}^{T}=2(m+1)(m-\frac{(1-t)}{4}\triangle_{\varphi_{t}}\varphi_{t}),
\end{eqnarray}
and
\begin{eqnarray}\label{MF4}
\begin{array}{lll}
\mathcal V_{d\eta} (\varphi_{t}  )
&=& -2(m+1)(I_{d\eta }-J_{d\eta})(\varphi_{t} ) +\frac{2}{V}\int_{M} h_{d\eta }d\eta ^{m}\wedge \eta  \\
&&-\frac{2t(m+1)}{V}\int_{M} \varphi_{t}(d\eta_{\varphi_{t}})^{m}\wedge \eta .
\end{array}
\end{eqnarray}
Then, by (\ref{2D}) and (\ref{Mabu}), we have
\begin{eqnarray}\label{MF5}
\begin{array}{lll}
\frac{d}{dt} \mathcal V_{d\eta }(\varphi_{t}) &=& -\frac{1}{V} \int_{M}
\dot{\varphi }_{t} (S ^{T}_{t}-\bar{S})(d\eta _{t})^{m}\wedge
\eta \\
&=& \frac{2m+2}{V} \int_{M}
\dot{\varphi }_{t} \frac{(1-t)}{4}\triangle_{\varphi_{t}}\varphi_{t}(d\eta _{t})^{m}\wedge
\eta \\
&=& 2(m+1)(t-1)\frac{d}{dt }((I_{d\eta }-J_{d\eta })(\varphi_{t})).
\end{array}
\end{eqnarray}
From (\ref{MF4}) and (\ref{MF5}), we have
\begin{eqnarray}\label{AA1}
\frac{d}{dt }(\frac{t}{V}\int_{M} \varphi_{t}(d\eta_{\varphi_{t}})^{m}\wedge \eta +t(I_{d\eta }-J_{d\eta })(\varphi_{t}))=(I_{d\eta }-J_{d\eta })(\varphi_{t}).
\end{eqnarray}
Noting that $\frac{t}{V}\int_{M} \varphi_{t}(d\eta_{\varphi_{t}})^{m}\wedge \eta +t(I_{d\eta }-J_{d\eta })(\varphi_{t})\rightarrow 0 $ as $t\rightarrow 0$. The  identity (\ref{AA1}) implies that
\begin{eqnarray}
\frac{1}{V}\int_{M} \varphi_{t}(d\eta_{\varphi_{t}})^{m}\wedge \eta +(I_{d\eta }-J_{d\eta })(\varphi_{t})\geq 0,
\end{eqnarray}
and
\begin{eqnarray}
\begin{array}{lll}
\mathcal V_{d\eta} (\varphi_{t}  )
&\leq & -2(m+1)(1-t)(I_{d\eta }-J_{d\eta})(\varphi_{t} ) +\frac{2}{V}\int_{M} h_{d\eta }d\eta ^{m}\wedge \eta  \\
&\leq & \frac{2}{V}\int_{M} h_{d\eta }d\eta ^{m}\wedge \eta .
\end{array}
\end{eqnarray}
Then the properness of $\mathcal V_{d\eta }$ implies that $J_{d\eta }(\varphi_{t})$, and consequently, $\|\varphi_{t}\|_{C^{0}}$ is uniformly bounded for $t\in [\epsilon , \tau )$. Therefore, the equation (\ref{MA1}) can also be solved.

\hfill $\Box$

\medskip

{\bf Proposition 4.4.} {\it Let  $(M, \xi , \eta , \Phi , g)$ be a compact Sasakian manifold with $[d\eta ]_{B} = \frac{2\pi }{m+1}
c_{1}^{B}(M, \mathcal F_{\xi})$. If $\mathcal V_{d\eta } $(or $F_{d\eta }$) is bounded from below in the space $\mathcal H(\xi , \eta , \Phi , g)$, then there exists a smooth family
of solution $\{\varphi_{t}\}$ of (\ref{MA2}) for $t\in (0, 1 )$.}

\medskip

{\bf Proof. } By corollary 3.5, it is sufficient to prove the $\mathcal K$-energy case. We prove it by contradiction. Let $\tau <1$ be the maximal number such that there exists a smooth family
of solution $\{\varphi_{t}\}$ of (\ref{MA2}) for $t\in (0, \tau )$. By (\ref{MF5}), we have
\begin{eqnarray}
\begin{array}{lll}
(I_{d\eta }-J_{d\eta })(\varphi_{t})&=&\int_{t_{0}}^{t}\frac{1}{2(m+1)(t-1)}\frac{d}{ds} \mathcal V_{d\eta }(\varphi_{s})+(I_{d\eta }-J_{d\eta })(\varphi_{t_{0}})\\
&\leq & \frac{1}{1-\tau } (V_{d\eta }(\varphi_{t_{0}})-\inf V_{d\eta }(\varphi))++(I_{d\eta }-J_{d\eta })(\varphi_{t_{0}}),
\end{array}
\end{eqnarray}
where $0<t_{0}<\tau $. Since $\mathcal V_{d\eta } $ is bounded from below, then  $(I_{d\eta }-J_{d\eta })(\varphi_{t})$ is bounded uniformly from above for $0<t<\tau$, and consequently, $\|\varphi_{t}\|_{C^{0}}$ is uniformly bounded for $t\in [t_{0} , \tau )$. So the equation (\ref{MA2}) can also be solved at $\tau $, this gives a contradiction.

\hfill $\Box$

\medskip

\section{Smoothing by the Sasakian-Ricci flow }
\setcounter{equation}{0}

\medskip

In this section, we use the Sasakian-Ricci flow to get a smoothing lemma. As a natural analogue of the K\"ahler-Ricci flow, the Sasakian-Ricci flow was introduced in \cite{SWZ}.
Now, we consider the following Sasakian-Ricci flow
\begin{eqnarray}\label{SR2}
\frac{\partial v}{\partial s }=\log \frac{(d\tilde{\eta}_{0}+\sqrt{-1}\partial_{B}\overline{\partial }_{B}v
)^{m}\wedge \tilde{\eta}_{0}}{(d\tilde{\eta}_{0} )^{m}\wedge \tilde{\eta }_{0}} +(m+1)v -h_{d\tilde{\eta}_{0}},
\end{eqnarray}
with $v|_{s=0}\equiv 0$. The long-time existence had been proved in \cite{SWZ}. In the following, for simplicity, we will denote the transverse K\"ahler form $d\tilde{\eta}_{0}+\sqrt{-1}\partial_{B}\overline{\partial }_{B}v$ by $d\tilde{\eta}_{s}$, and we will use a subscript $s$ to indicated objects that are defined with respect to the transverse K\"ahler metric $d\tilde{\eta}_{s}$.  As that in \cite{Bando}, we have the following lemma.

\medskip

{\bf Lemma 5.1. } {\it The following inequalities
\begin{eqnarray}\label{SRL1}
\|\frac{\partial v}{\partial s}\|_{C^{}0}\leq e^{(m+1)s}\|h_{d\tilde{\eta}_{0}}\|_{C^{0}},
\end{eqnarray}
\begin{eqnarray}\label{SRL2}
\sup_{M}(|h_{d\tilde{\eta}_{s}}|^{2}+\frac{s}{2}|d h_{d\tilde{\eta}_{s}}|_{s}^{2})\leq 4e^{2(m+1)s}\|h_{d\tilde{\eta}_{0}}\|_{C^{0}}^{2},
\end{eqnarray}
\begin{eqnarray}\label{SRL3}
e^{-(m+1)s}\triangle_{s}h_{d\tilde{\eta}_{s}}\geq \triangle_{0}h_{d\tilde{\eta}_{0}},
\end{eqnarray}
hold for all $s\geq 0$.

}

\medskip

{\bf Proof. } Differentiating  the Sasakian-Ricci flow equation (\ref{SR2}) gives
\begin{eqnarray}\label{SRR1}
\frac{\partial }{\partial s} \dot{v}=\frac{1}{4}\triangle_{s} \dot{v}+(m+1)\dot{v},
\end{eqnarray}
 the maximum principle implies (\ref{SRL1}).  By direct calculation, we have
\begin{eqnarray}\label{SRR2}
\frac{\partial }{\partial s} |d \dot{v}|_{s}^{2}=\frac{1}{4}\triangle_{s}|d \dot{v}|_{s}^{2} -\frac{1}{2}|\nabla_{s}d \dot{v}|_{s}^{2}+(m+1)|d \dot{v}|_{s}^{2}.
\end{eqnarray}
From (\ref{SRR1}) and (\ref{SRR2}), we have
 \begin{eqnarray}\label{SRR3}
(\frac{\partial }{\partial s}-\frac{1}{4}\triangle_{s} )(\dot{v}^{2}+\frac{s}{2}|d \dot{v}|_{s}^{2})\leq 2(m+1)(\dot{v}^{2}+\frac{s}{2}|d \dot{v}|_{s}^{2}),
\end{eqnarray}
and the maximum principle implies that
\begin{eqnarray}\label{SRR4}
\sup_{M}(\dot{v}^{2}+\frac{s}{2}|d \dot{v}|_{s}^{2})\leq e^{2(m+1)s}\|h_{d\tilde{\eta}_{0}}\|_{C^{0}}^{2}.
\end{eqnarray}
By the equation (\ref{SR2}), it is easy to check that
\begin{eqnarray}
h_{d\tilde{\eta}_{s}} =-\dot{v} +c_{s}
\end{eqnarray}
for some constant $c_{s}$ with $c_{0}=0$. From the normalization condition the Ricci potential function and (\ref{SR2}), we have
\begin{eqnarray}
\int_{M} e^{-(m+1)v +h_{d\tilde{\eta}_{0}}+c_{s}} (d\tilde{\eta}_{0})^{m}\wedge \tilde{\eta}_{0}=\int_{M} e^{h_{d\tilde{\eta}_{s}}} (d\tilde{\eta}_{s})^{m}\wedge \tilde{\eta}_{s}=V,
\end{eqnarray}
and then
\begin{eqnarray}\label{SRR5}
|c_{s}|\leq (m+1)\|v\|_{C^{0}}\leq e^{(m+1)s}\|h_{d\tilde{\eta}_{0}}\|_{C^{0}}.
\end{eqnarray}
So, (\ref{SRR4}) and (\ref{SRR5}) imply (\ref{SRL2}).

By direct calculation, one can check that
\begin{eqnarray}\label{SRR6}
(\frac{\partial }{\partial s}-\frac{1}{4}\triangle_{s} )(\triangle_{s}\dot{v}) =(m+1)\triangle_{s}\dot{v}-|\partial_{B}\bar{\partial}_{B}\dot{v}|_{s}^{2},
\end{eqnarray}
and then the maximum principle implies the inequality (\ref{SRL3}).

\hfill $\Box$

\medskip

{\bf Lemma 5.2. } {\it Let $v_{t, s}$ be a solution of (\ref{SR2}) with $d\tilde{\eta}_{0}= d\eta_{\varphi _{t}}$. Let $\tilde{h}=h_{d\tilde{\eta}_{1}} -\frac{1}{V}\int_{M} h_{d\tilde{\eta}_{1}} (d\tilde{\eta}_{1})^{m}\wedge \tilde{\eta}_{1} $ and assume that
\begin{eqnarray}\label{V1}
\frac{1}{2}d\eta_{SE}\leq d\tilde{\eta}_{1}\leq d\eta_{SE}.
\end{eqnarray}
Then for any $p> 2m+1$, there exist positive constant $\bar{C}_{1}$ depending only on $(M, g_{SE})$ and $p$ such that
\begin{eqnarray}\label{ZZZ3}
\|\tilde{h}\|_{C^{0}}\leq \bar{C}_{1}(1-t)^{\frac{1}{p-1}}\|h_{d\tilde{\eta}_{0}}\|_{C^{0}}^{\frac{p-2}{p-1}}.
\end{eqnarray}
}

\medskip

{\bf Proof. } (\ref{SRL2}) implies that $\|\tilde{h}\|_{C^{0}}\leq 4e^{m+1} \|h_{d\tilde{\eta}_{0}}\|_{C^{0}} $. By the initial condition $d\tilde{\eta}_{0}= d\eta_{\varphi _{t}}$, we have $\rho^{T}(d\tilde{\eta}_{0})\geq t (m+1)d\tilde{\eta}_{0}$, and $\triangle_{0} h_{d\tilde{\eta}_{0}}\geq 4 m(m+1)(t-1)$. By (\ref{SRL3}), we have
\begin{eqnarray}
-\triangle_{1}h_{d\tilde{\eta}_{1}}\leq 4e^{m+1}m(m+1)(1-t).
\end{eqnarray}
Integrating by parts, we have
\begin{eqnarray}\label{ZZZ1}
\begin{array}{lll}
\int_{M}|d \tilde{h}|_{1}^{2}(d\tilde{\eta}_{1})^{m}\wedge \tilde{\eta}_{1} &= & -\int_{M}\tilde{h}\triangle_{1}\tilde{h}(d\tilde{\eta}_{1})^{m}\wedge \tilde{\eta}_{1}\\
&\leq & \int_{M}(\tilde{h} - \inf \tilde{h})\sup_{M}(-\triangle_{1}\tilde{h})(d\tilde{\eta}_{1})^{m}\wedge \tilde{\eta}_{1}\\
&\leq & 2V \|\tilde{h}\|_{C^{0}}\sup_{M}(-\triangle_{1}\tilde{h})\\
&\leq & \bar{C}_{2}(1-t)\|\tilde{h}\|_{C^{0}},\\
\end{array}
\end{eqnarray}
where $\bar{C}_{2}$ depends only on the dimension of $M$. Since $\tilde{h}$ be a basic function, by condition (\ref{V1}), we have
\begin{eqnarray}
|d\tilde{h}|_{SE}^{2}\leq 2 |d\tilde{h}|_{1}^{2}
\end{eqnarray}
 Let $p> 2m+1$, by the Sobolev imbedding theorem (Lemma 2.22 of \cite{Au2}), the Poincar\'e inequality and (\ref{SRL2}), we have
\begin{eqnarray}\label{ZZZ2}
\begin{array}{lll}
\|\tilde{h}\|_{C^{0}}^{p}&\leq & \bar{C}_{3}(\int_{M}|\tilde{h}|^{p}+|d \tilde{h}|_{SE}^{p}(d\eta_{SE})^{m}\wedge \eta_{SE})\\
&\leq & \bar{C}_{4}\|h_{d\tilde{\eta}_{0}}\|_{C^{0}}^{p-2}(\int_{M}|\tilde{h}|^{2}+|d \tilde{h}|_{SE}^{2}(d\eta_{SE})^{m}\wedge \eta_{SE})\\
&\leq & \bar{C}_{5}\|h_{d\tilde{\eta}_{0}}\|_{C^{0}}^{p-2}(\int_{M}|d \tilde{h}|_{SE}^{2}(d\eta_{SE})^{m}\wedge \eta_{SE})\\
&\leq & \bar{C}_{6}\|h_{d\tilde{\eta}_{0}}\|_{C^{0}}^{p-2}(\int_{M}|d \tilde{h}|_{1}^{2}(d\tilde{\eta }_{1})^{m}\wedge \tilde{\eta}_{1}),\\
\end{array}
\end{eqnarray}
where  constants $\bar{C}_{i}$ depends only on $(M, g_{SE})$ and $p$. Then (\ref{ZZZ1}) and (\ref{ZZZ2}) imply  (\ref{ZZZ3}), and we are finished.

\hfill $\Box$

\medskip

{\bf Lemma 5.3. } {\it Let $v_{t, s}$ be a solution of (\ref{SR2}) with initial data $d\tilde{\eta}_{0}= d\eta_{\varphi _{t}}$, and $u_{t}=v_{t, 1}$. We have the inequality
\begin{eqnarray}\label{531}
\|u_{t}\|_{C^{0}} \leq \frac{1}{m+1}e^{m+1}\|h_{d\eta_{\varphi_{t}}}\|_{C^{0}}
\end{eqnarray}
for all $t\in [0, 1]$. Moreover, assume that $\frac{1}{2}d\eta_{SE}\leq d\eta_{\varphi_{t}+u_{t}}\leq d\eta_{SE}$ for all $t\in [t_{1} , 1]$, where $t_{1}\in [0, 1)$. Then for any $p>2m+1$ and $0\leq k<1$, there exists a constant $\bar{C}_{7}$ depending only on $(M, g_{SE})$ and $p$ such that
\begin{eqnarray}\label{532}
\|h_{d\eta_{\varphi_{t}+u_{t}}}\|_{C^{0,k}(d\eta_{SE})}\leq \bar{C}_{7}(1-t)^{1-\beta  }(1+\|h_{d\eta_{\varphi_{t}}}\|_{C^{0}})^{\beta}
\end{eqnarray}
for all $t\in [t_{1} , 1]$, where $\beta =\frac{p+k-2}{p-1}$.

}

\medskip

{\bf Proof. } From (\ref{SRL1}), it follows that $|\frac{\partial v_{t, s}}{\partial s}|\leq e^{(m+1)s}\|h_{d\eta_{\varphi_{t}}}\|_{C^{0}}$, and integrating from $0$ to $1$, we obtain the inequality (\ref{531}).

In the following, let $d (x, y)$ be the distance between $x$ and $y$ with respect to the metric $g_{SE}$. Since $h_{d\eta_{\varphi_{t}+u_{t}}}$ is a basic function, by the condition $\frac{1}{2}d\eta_{SE}\leq d\eta_{\varphi_{t}+u_{t}}\leq d\eta_{SE}$, we have $$|dh_{d\eta_{\varphi_{t}+u_{t}}}|_{d\eta_{SE}}\leq \sqrt{2} |dh_{d\eta_{\varphi_{t}+u_{t}}}|_{d\eta_{\varphi_{t}+u_{t}}}.$$
If $d(x, y )\leq (1-t)^{\frac{1}{p-1}}(1+\|h_{d\eta_{\varphi_{t}}}\|_{C^{0}})^{-\frac{1}{p-1}}$, by (\ref{SRL2}) in lemma 5.1, we have
\begin{eqnarray}\label{533}
\begin{array}{lll}
&&|h_{d\eta_{\varphi_{t}+u_{t}}}(x)-h_{d\eta_{\varphi_{t}+u_{t}}}(y)| \leq d(x, y)\sup_{M}|dh_{d\eta_{\varphi_{t}+u_{t}}}|_{d\eta_{SE}}\\
&\leq & \sqrt{2}d(x, y)\sup_{M}|dh_{d\eta_{\varphi_{t}+u_{t}}}|_{d\eta_{\varphi_{t}+u_{t}}}\\
&\leq & 4\sqrt{2}e^{m+1} d(x, y)(1+\|h_{d\eta_{\varphi_{t}}}\|_{C^{0}})\\
&\leq & 4\sqrt{2}e^{m+1} (1-t)^{\frac{1-k}{p-1}}(1+\|h_{d\eta_{\varphi_{t}}}\|_{C^{0}})^{\frac{p+k-2}{p-1}}d(x, y)^{k}.\\
\end{array}
\end{eqnarray}
If $d(x, y )\geq (1-t)^{\frac{1}{p-1}}(1+\|h_{d\eta_{\varphi_{t}}}\|_{C^{0}})^{-\frac{1}{p-1}}$, then the estimate (\ref{ZZZ3}) in lemma 5.2 implies
\begin{eqnarray}\label{534}
\begin{array}{lll}
&&|h_{d\eta_{\varphi_{t}+u_{t}}}(x)-h_{d\eta_{\varphi_{t}+u_{t}}}(y)| \leq 2\|\tilde{h}\|_{C^{0}}\\
&\leq & 2\bar{C}_{1}(1-t)^{\frac{1}{p-1}}(\|h_{d\eta_{\varphi_{t}}}\|_{C^{0}})^{\frac{p-2}{p-1}}\\
&\leq & 2\bar{C}_{1} (1-t)^{\frac{1-k}{p-1}}(1+\|h_{d\eta_{\varphi_{t}}}\|_{C^{0}})^{\frac{p+k-2}{p-1}}d(x, y)^{k}.\\
\end{array}
\end{eqnarray}
On the other hand, the integral  normalization $\int_{M}e^{h_{d\eta_{\varphi_{t}+u_{t}}}} (d\eta_{\varphi_{t}+u_{t}})^{m}\wedge \eta =V$ implies $h_{d\eta_{\varphi_{t}+u_{t}}}$ change signs, so we  have
\begin{eqnarray}\label{535}
\begin{array}{lll}
&&\|h_{d\eta_{\varphi_{t}+u_{t}}}\|_{C^{0}} \leq Osc (h_{d\eta_{\varphi_{t}+u_{t}}})= Osc (\tilde{h})\leq 2\|\tilde{h}\|_{C^{0}}\\
&\leq & 2\bar{C}_{1}(1-t)^{\frac{1}{p-1}}(\|h_{d\eta_{\varphi_{t}}}\|_{C^{0}})^{\frac{p-2}{p-1}}.\\
\end{array}
\end{eqnarray}
It is easy to see that (\ref{533}), (\ref{534}) and (\ref{535}) imply the estimate (\ref{532}).

\hfill $\Box$

\medskip

Set $\alpha :=1-\frac{1}{4m+2}>\frac{1}{2}$ and define the function $f_{d\eta }$ by
\begin{eqnarray}\label{DD}
f_{d\eta }(t):=(1-t)^{1-\alpha }(1+2(1-t)\|\varphi_{t}\|_{C^{0}})^{\alpha }.\end{eqnarray}

Discussing as that in \cite{T1}, we have the following proposition.

\medskip

{\bf Proposition 5.4. } {\it Suppose that $(M, \xi, \eta, \Phi, g)$ admits no non-trivial Hamiltonian holomorphic vector fields. Let $\varphi _{t}$ be a smooth family of solutions of the equation (\ref{MA2}) for $t\in (0, 1]$. There exist a constant $D>0$ depending only on $(M, g_{SE})$ such that
\begin{eqnarray}\label{PP}
\|\varphi_{1}- \varphi_{t}\|_{C^{0}}\leq A (1-t) \|\varphi_{t}\|_{C^{0}}+1
\end{eqnarray}
for all $t\in [t_{0} , 1]$, where $t_{0}\in [0, 1)$ satisfies $f_{d\eta } (t_{0})=\max _{[t_{0} , 1]}f_{d\eta } =D$ and $A$ depending only on the dimension of $M$.

}

\medskip

{\bf Proof. } Let's rewrite (\ref{MA2}) as the following transverse Monge-Amp\`ere equation with $d\eta_{SE}$ as reference metric
\begin{eqnarray}\label{MA31}
\begin{array}{lll}
&&\frac{(d\eta_{SE}+\sqrt{-1}\partial_{B}\overline{\partial }_{B}(\varphi_{t}-\varphi_{1} )
)^{m}\wedge \eta}{(d\eta_{SE} )^{m}\wedge \eta}\\&&=\exp (-(m+1)(\varphi_{t}-\varphi_{1} )+(1-t)(m+1)\varphi_{t}
).\\
\end{array}
\end{eqnarray}
 It is easy to see that $h_{d\eta _{\varphi_{t}}}=(t-1)(m+1)\varphi_{t}+c_{t}$, for some constant $c_{t}$. The integrate normalization of the Ricci potential function $h_{d\eta_{\varphi_{t}}}$ gives
\begin{eqnarray}
\begin{array}{lll}
V&=&\int_{M}(d\eta_{\varphi_{t}})^{m}\wedge \eta =\int_{M}e^{h_{d\eta_{\varphi_{t}}}}(d\eta_{\varphi_{t}})^{m}\wedge \eta \\
&=&\int_{M}e^{(t-1)(m+1)\varphi_{t}+c_{t}}(d\eta_{\varphi_{t}})^{m}\wedge \eta ,\\
\end{array}
\end{eqnarray}
from which it follows that
\begin{eqnarray}\label{541}
|c_{t}|\leq (m+1)(1-t)\|\varphi_{t}\|_{C^{0}},
\end{eqnarray}
and
\begin{eqnarray}\label{545}
\|h_{d\eta_{\varphi_{t}}}\|_{C^{0}}\leq 2(m+1)(1-t)\|\varphi_{t}\|_{C^{0}}.
\end{eqnarray}
Then, lemma 5.3 implies that
\begin{eqnarray}\label{542}
\|u_{t}\|_{C^{0}}\leq 2e^{(m+1)}(1-t)\|\varphi_{t}\|_{C^{0}}.
\end{eqnarray}

Consider $d\eta_{\varphi_{t}+u_{t}}=d\eta +\sqrt{-1}\partial_{B}\bar{\partial }_{B}(\varphi_{t}+u_{t})=d\eta_{SE}+\sqrt{-1}\partial_{B}\bar{\partial }_{B}(\varphi_{t}+u_{t}-\varphi_{1})$, and then
\begin{eqnarray}\label{MA312}
\begin{array}{lll}
&&\frac{(d\eta_{SE}+\sqrt{-1}\partial_{B}\overline{\partial }_{B}(\varphi_{t}+u_{t}-\varphi_{1} )
)^{m}\wedge \eta}{(d\eta_{SE} )^{m}\wedge \eta}\\&&=\exp (-(m+1)(\varphi_{t}+u_{t}-\varphi_{1} )-h_{d\eta_{\varphi_{t}+u_{t}}}-\tilde{c}_{t}
)\\
\end{array}
\end{eqnarray}
for some constant $\tilde{c}_{t}$.  Setting $\tilde{\varphi}_{t}=\varphi_{t}+u_{t}-\varphi_{1}+\frac{\tilde{c}_{t}}{m+1}$, from (\ref{MA312}) and (\ref{542}), we have \begin{eqnarray}
\begin{array}{lll}
&&\int_{M}e^{h_{d\eta_{\varphi_{t}+u_{t}}}}(d\eta_{\varphi_{t}+u_{t}})^{m}\wedge \eta =\int_{M}e^{-(m+1)\tilde{\varphi}_{t}}(d\eta_{SE})^{m}\wedge \eta \\
&=&\int_{M}e^{-(m+1)\tilde{\varphi}_{t}+t(m+1)\varphi_{t}-(m+1)\varphi_{1}}(d\eta_{\varphi_{t}})^{m}\wedge \eta \\
&=&\int_{M}e^{(t-1)(m+1)\varphi_{t}-(m+1)u_{t}-\tilde{c}_{t}}(d\eta_{\varphi_{t}})^{m}\wedge \eta ,\\
\end{array}
\end{eqnarray}
and then
\begin{eqnarray}\label{543}
\begin{array}{lll}
|\tilde{c}_{t}|&\leq & (1-t)(m+1)\|\varphi_{t}\|_{C^{0}}+(m+1)\|u_{t}\|_{C^{0}}\\
&\leq & (1-t)(m+1)(1+2e^{(m+1)})\|\varphi_{t}\|_{C^{0}}.
\end{array}
\end{eqnarray}

Recall that $\varphi_{t}-\varphi_{1}=\tilde{\varphi}_{t}-u_{t}-\frac{\tilde{c}_{t}}{m+1}$, from (\ref{542}) and (\ref{543}), we have
\begin{eqnarray}
\|\varphi_{t}-\varphi_{1}\|_{C^{0}}=\|\tilde{\varphi}_{t}\|_{C^{0}}+(1-t)(4e^{(m+1)}+1)\|\varphi_{t}\|_{C^{0}}.
\end{eqnarray}
From above, it will suffice to get the estimate $\|\tilde{\varphi}_{t}\|_{C^{0}}\leq 1$.

Let's consider the following transverse Monge-Amp\`ere equation
\begin{eqnarray}\label{MA4}
\log \{\frac{(d\eta_{SE}+\sqrt{-1}\partial_{B}\overline{\partial }_{B}\psi
)^{m}\wedge \eta}{(d\eta_{SE} )^{m}\wedge \eta}\}+(m+1)\psi =\tilde{\psi} .
\end{eqnarray}
The linearization of the left side of (\ref{MA4}) at $\psi =0$ is
\begin{eqnarray}
\delta \psi \mapsto \frac{1}{4}\triangle_{SE} \delta \psi + (m+1)\delta \psi ,
\end{eqnarray}
which is a transverse elliptic operator from $C_{B}^{i+2 , k } (M)\rightarrow C_{B}^{i+2 , k } (M)$ for any $0<k <1$ and $i\geq 0$. If $M$ doesn't have non-trivial Hamiltonian holomorphic vector fields, by theorem 5.1 of \cite{12}, we have $ker (\frac{1}{4}\triangle_{SE} + (m+1))=0$, then the operator $(\frac{1}{4}\triangle_{SE} + (m+1)): C_{B}^{i+2 , \epsilon } (M)\rightarrow C_{B}^{i+2 , \epsilon } (M)$ is invertible. Applying the implicit function theorem, there exist positive constants $\epsilon(d\eta_{SE})$ and $C^{\ast}(d\eta_{SE})$ which depend only on $k$ and the geometry of $(M, g_{SE})$, so that
\begin{eqnarray}\label{544}
if \quad \|\tilde{\psi } \|_{C^{0, k }}\leq \epsilon(d\eta_{SE}) \quad  then \quad \|\psi\|_{C^{2, k }}\leq C^{\ast }(d\eta_{SE})\|\tilde{\psi}\|_{C^{0, k }}.\end{eqnarray}

Setting $D =\frac{\epsilon (m+1)^{-\alpha}}{2(\bar{C}_{7}+1)(C^{\ast}+1)(\epsilon +1)}$, where $\epsilon =\epsilon(d\eta_{SE})$, $C^{\ast}=C^{\ast }(d\eta_{SE})$ are chosen as in (\ref{544}), $\alpha =1-\frac{1}{4m+2}$, $\bar{C}_{7}$ is defined as in lemma 5.3 (by choosing $k=\frac{1}{2}$ and $p=2m+2$. Let $t_{0}\in [0, 1)$ satisfies $f_{d\eta } (t_{0})=\max _{[t_{0} , 1]}f_{d\eta } =D$. Now, we only need to prove the following claim.

\medskip

{\bf Claim } {For all $t\in [t_{0} , 1]$, we have
\begin{eqnarray}
\|\tilde{\varphi}_{t}\|_{C^{2, \frac{1}{2} }}<\frac{1}{2}.
\end{eqnarray}

}

\medskip

We assume the contrary. Since $\tilde{\varphi_{1}}=0$, there exists $t_{1}\in [t_{0}, 1)$ such that
\begin{eqnarray}
\|\tilde{\varphi}_{t_1}\|_{C^{2 , \frac{1}{2}}(d\eta_{SE})}=\frac{1}{2}, \quad and \quad \|\tilde{\varphi}_{t}\|_{C^{2 , \frac{1}{2}}(d\eta_{SE})}<\frac{1}{2} \quad if \quad t_{1}<t<1.
\end{eqnarray}
In particular $-\frac{1}{4}d\eta_{SE}\leq \sqrt{-1}\partial_{B}\overline{\partial }_{B}\tilde{\varphi}_{t}\leq \frac{1}{4}d\eta_{SE}$, and then
\begin{eqnarray}
\frac{3}{4}d\eta_{SE}\leq d\eta_{\varphi_{t}+u_{t}}\leq \frac{5}{4}d\eta_{SE}
\end{eqnarray}
for all $t\in [t_{1}, 1]$. By applying (\ref{532}) in lemma 5.3 (by choosing $p=2m+2$) and (\ref{545}), we have
\begin{eqnarray}
\begin{array}{lll}
\|h_{d\eta_{\varphi_{t}+u_{t}}}\|_{C^{0,\frac{1}{2}}(d\eta_{SE})}&\leq & \bar{C}_{7}(1-t)^{1-\alpha  }(1+\|h_{d\eta_{\varphi_{t}}}\|_{C^{0}})^{\alpha}\\
&\leq & \bar{C}_{7}(1-t)^{1-\alpha  }(1+2(1-t)(m+1)\|\varphi_{t}\|_{C^{0}})^{\alpha}\\
&\leq & \bar{C}_{7}(m+1)^{\alpha }(1-t)^{1-\alpha  }(1+2(1-t)\|\varphi_{t}\|_{C^{0}})^{\alpha}\\
&\leq & \bar{C}_{7}(m+1)^{\alpha }D\\
&=& \frac{\bar{C}_{7}\epsilon }{2(\bar{C}_{7}+1)(C^{\ast}+1)(\epsilon +1)}\\
&<& \epsilon ,
\end{array}
\end{eqnarray}
for all $t\in [t_{1}, 1]$. Using (\ref{544}) again, we get
\begin{eqnarray}
\begin{array}{lll}
\|\tilde{\varphi}_{t_1}\|_{C^{2 , \frac{1}{2}}(d\eta_{SE})}&\leq & \|h_{d\eta_{\varphi_{t}+u_{t}}}\|_{C^{0,\frac{1}{2}}(d\eta_{SE})}\\
&\leq & \frac{C^{\ast}\bar{C}_{7}\epsilon }{2(\bar{C}_{7}+1)(C^{\ast}+1)(\epsilon +1)}\\
&<& \frac{1}{2}.
\end{array}
\end{eqnarray}
This gives a contradiction, and complete the proof of the claim. So, the proof of the proposition is complete.

\hfill $\Box$

\section{A Moser-Trudinger type inequality }
\setcounter{equation}{0}

In this section, we assume the existence of a Sasakian-Einstein structure  and establish a Moser-Trudinger type inequality for functional $F_{d\eta_{SE}}$, our discussion follow that in \cite{PSSW} by Phong, Song, Strum and Weinkove. In fact, we obtain the following theorem.

 \medskip

{\bf Theorem 6.1. } {\it
Let $(M, \xi , \eta_{SE}, \Phi_{SE}, g_{SE})$ be a compact Sasakian-Einstein metric  without non trivial Hamiltonian holomorphic vector field, then there exist uniform positive constants $C_{1}$, $C_{2}$ depending only the geometry of $(M, g_{SE})$, such that
\begin{eqnarray}\label{MT2}
F_{d\eta_{SE} }(\varphi) \geq C_{1}J_{d\eta_{SE} }(\varphi ) -C_{2},
\end{eqnarray}
for all $\varphi \in \mathcal H (\xi , \eta_{SE} , \Phi_{SE} , g_{SE})$.}

\medskip

{\bf Proof. }
 Fix a basic function $\phi \in \mathcal H (\xi , \eta_{E} , \Phi_{SE} , g_{SE})$, and set $d\eta =d\eta_{SE}+\sqrt{-1}\partial_{B}\bar{\partial }_{B}\phi $. Now, let us consider the complex Monge-Amp\`ere equation (\ref{MA2}). Since there are no nontrivial Hamiltonian holomorphic vector fields, by the uniqueness of Sasakian-Einstein structure (\cite{Se} or \cite{Ni}) and  proposition 4.1, a unique solution $\varphi_{t}$ exists for all $t\in (0, 1]$, and $d\eta_{\varphi_{1}}=d\eta_{SE}$. In particular $\varphi_{1}$ and $-\phi $ differ by a constant.

For further consideration, we give the following estimates for functionals $F$, $I$ and $J$. From (\ref{DF}), (\ref{DJ}) and (\ref{DM}), we have
\begin{eqnarray}
\begin{array}{lll}
\frac{d}{ds }(I_{d\eta }-J_{d\eta })(\varphi_{s})&=& -\frac{d}{ds }(\frac{1}{V}\int_{M}\varphi_{s}(d\eta_{\varphi_{s}})^{m}\wedge \eta )
+\frac{1}{V}\int_{M}\dot{\varphi }_{s}(d\eta)^{m}\wedge \eta
\\
&& -\frac{1}{V}\int_{M}\dot{\varphi }_{s}\{(d\eta)^{m}-(d\eta_{\varphi_{s}})^{m}\}\wedge \eta \\
&=& -\frac{d}{ds }(\frac{1}{V}\int_{M}\varphi_{s}(d\eta_{\varphi_{s}})^{m}\wedge \eta )
+\frac{1}{V}\int_{M}\dot{\varphi }_{s}(d\eta_{\varphi_{s}})^{m}\wedge \eta
\\
&=& -\frac{d}{ds }(\frac{1}{V}\int_{M}\varphi_{s}(d\eta_{\varphi_{s}})^{m}\wedge \eta )
-\frac{1}{sV}\int_{M}\varphi_{s} (d\eta_{\varphi_{s}})^{m}\wedge \eta .
\\
\end{array}
\end{eqnarray}
The uniform $C^{0}$ estimate (\ref{C02}) of $\varphi_{t}$ implies that
\begin{eqnarray}
s\frac{1}{V}\int_{M}\varphi_{s}(d\eta_{\varphi_{s}})^{m}\wedge \eta \rightarrow 0
\end{eqnarray}
as $s\rightarrow 0$.
By integrating on $[0, t]$, we get
\begin{eqnarray}\label{K1}
\begin{array}{lll}
&&t(I_{d\eta }-J_{d\eta })(\varphi_{t})-\int_{0}^{t}(I_{d\eta }-J_{d\eta })(\varphi_{s}) ds\\
&=&
\int_{0}^{t}s\frac{d}{ds }(I_{d\eta }-J_{d\eta })(\varphi_{s})ds\\
&=& -\int_{0}^{t}s\frac{d}{ds }(\frac{1}{V}\int_{M}\varphi_{s}(d\eta_{\varphi_{s}})^{m}\wedge \eta )ds
-\frac{1}{V}\int_{0}^{t}(\int_{M}\varphi_{s} (d\eta_{\varphi_{s}})^{m}\wedge \eta )ds
\\
&=&-\frac{t}{V}\int_{M}\varphi_{t} (d\eta_{\varphi_{t}})^{m}\wedge \eta ,\\
\end{array}
\end{eqnarray}
and then
\begin{eqnarray}\label{K2}
\begin{array}{lll}
F_{d\eta }^{0}(\varphi_{t})&=&-(I_{d\eta }-J_{d\eta })(\varphi_{t})
-\frac{1}{V}\int_{M}\varphi_{t} (d\eta_{\varphi_{t}})^{m}\wedge \eta\\
&=& \frac{-1}{t}\int_{0}^{t}(I_{d\eta }-J_{d\eta })(\varphi_{s}) ds .\\
\end{array}
\end{eqnarray}
Taking $t=1$ and considering $F_{d\eta }(\varphi_{1})=-F_{d\eta_{SE}}(\phi )$, so that
\begin{eqnarray}\label{SE1}
F_{d\eta_{SE} }(\phi)=\int_{0}^{1}(I_{d\eta }-J_{d\eta })(\varphi_{s}) ds .
\end{eqnarray}

By the definitions (\ref{DF}) and the cocycle property of $F_{d\eta }^{0}$, we have
\begin{eqnarray}\label{J1}
\begin{array}{lll}
J_{d\eta }(\varphi_{1})-J_{d\eta }(\varphi_{t})&=& \frac{1}{V}\int_{M}(\varphi_{1}-\varphi_{t})(d\eta)^{m}\wedge \eta
+F_{d\eta }^{0}(\varphi_{1}) -F_{d\eta }^{0}(\varphi_{t})\\
&=& \frac{1}{V}\int_{M}(\varphi_{1}-\varphi_{t})(d\eta)^{m}\wedge \eta
-F_{d\eta _{\varphi_{1}}}^{0}(\varphi_{t}-\varphi_{1})\\
&\leq & \frac{1}{V}\int_{M}(\varphi_{1}-\varphi_{t})(d\eta)^{m}\wedge \eta
+\frac{1}{V}\int_{M}(\varphi_{t}-\varphi_{1})(d\eta_{\varphi_{1}})^{m}\wedge \eta)\\
&\leq & Osc (\varphi_{1}-\varphi_{t}).
\end{array}
\end{eqnarray}
By adding
\begin{eqnarray*}
\begin{array}{lll}
&&I_{d\eta }(\varphi_{t})-I_{d\eta }(\varphi_{1})\\&=& \frac{1}{V}\int_{M}(\varphi_{t}-\varphi_{1})(d\eta)^{m}\wedge \eta + \frac{1}{V}\int_{M}(\varphi_{1}-\varphi_{t})(d\eta_{\varphi_{1}}^{m})\wedge \eta
\\
& &
+\frac{1}{V}\int_{M}\varphi_{t}\{(d\eta_{\varphi_{1}})^{m}-(d\eta_{\varphi_{t}})^{m}\}\wedge \eta ,\\
\end{array}
\end{eqnarray*}
we get
\begin{eqnarray}\label{J2}
\begin{array}{lll}
&&(I_{d\eta } -J_{d\eta })(\varphi_{t})-(I_{d\eta} -J_{d\eta })(\varphi_{1})\\
&=&
J_{d\eta }(\varphi_{1})-J_{d\eta }(\varphi_{t})+(I_{d\eta }(\varphi_{t})-I_{d\eta }(\varphi_{1}))\\
&\leq & \frac{1}{V}\int_{M}\varphi_{t}\{(\omega_{\varphi_{1}})^{m}-(\omega_{\varphi_{t}})^{m}\} \\
&=& \frac{1}{V}\int_{M}\varphi_{t}(d\eta_{\varphi_{1}}-d\eta_{\varphi_{t}})\wedge (\sum_{j=0}^{m-1}d\eta_{\varphi_{t}}^{j}\wedge d\eta_{\varphi_{1}}^{(m-1-j)})\wedge \eta \\
&=& \frac{1}{V}\int_{M}(\varphi_{1}-\varphi_{t})(d\eta_{\varphi_{t}}-d\eta)\wedge (\sum_{j=0}^{m-1}d\eta_{\varphi_{t}}^{j}\wedge d\eta_{\varphi_{1}}^{(m-1-j)})\wedge \eta \\
&\leq & m Osc (\varphi_{1}-\varphi_{t}).\\
\end{array}
\end{eqnarray}
Interchanging $\varphi_{t}$ and $\varphi_{1}$ in (\ref{J1}) and (\ref{J2}), we get
\begin{eqnarray}\label{J3}
|J_{d\eta }(\varphi_{1})-J_{d\eta }(\varphi_{t})|\leq Osc (\varphi_{1}-\varphi_{t})
\end{eqnarray}
and
\begin{eqnarray}\label{J4}
|(I_{d\eta } -J_{d\eta})(\varphi_{t})-(I_{d\eta} -J_{d\eta })(\varphi_{1})|\leq m \cdot Osc (\varphi_{1}-\varphi_{t}).
\end{eqnarray}

Using the relationship $F_{d\eta }(\varphi_{1})=-F_{d\eta_{SE}}(\phi)$, we have
\begin{eqnarray}\label{J5}
\begin{array}{lll}
J_{d\eta}(\varphi_{1}) &=& F_{d\eta }(\varphi_{1}) +\frac{1}{V}\int_{M}\varphi_{1}(d\eta
)^{m} \wedge \eta \\
&=& -F_{d\eta_{SE}}(\phi)+\frac{1}{V}\int_{M}\varphi_{1}(d\eta
)^{m}\wedge \eta \\
&=&-J_{d\eta_{SE}}(\phi)+\frac{1}{V}\int_{M}\phi\{(d\eta_{SE}
)^{m}-(d\eta)^{m}\wedge \eta \} \\
&=&(I_{d\eta_{SE}}-J_{d\eta_{SE}})(\phi)\geq \frac{1}{m} J_{d\eta_{SE}}(\phi),
\end{array}
\end{eqnarray}
where we have used the inequality (\ref{2.17}).
Since $(I_{d\eta }-J_{d\eta})(\varphi_{t})$ is nondecreasing in $t$, (\ref{SE1}) implies that
\begin{eqnarray}
F_{d\eta_{SE} }(\phi)\geq (1-t)(I_{d\eta }-J_{d\eta })(\varphi_{t}) ds \geq \frac{1-t}{m}J_{d\eta}(\varphi_{t}),
\end{eqnarray}
using (\ref{J5}) and (\ref{J1}),  we have
\begin{eqnarray}\label{F3}
F_{d\eta_{SE} }(\phi)\geq  \frac{1-t}{m^{2}}J_{d\eta_{SE} }(\phi )-\frac{1-t}{m}Osc (\varphi_{t}-\varphi_{1}).
\end{eqnarray}

\medskip

In the following, we choose $t_{0}$ as that in proposition 5.4. If $2(1-t_{0})\|\varphi_{t_{0}}\|_{C^{0}}\leq 1$, by the definition of $t_{0}$ , we have $D\leq (1-t_{0})^{1-\alpha }2^{\alpha}$, i.e.
\begin{eqnarray}
(1-t_{0})\geq 2^{-\frac{\alpha }{1-\alpha }}D^{\frac{1 }{1-\alpha }}.
\end{eqnarray}
If $2(1-t_{0})\|\varphi_{t_{0}}\|_{C^{0}}\geq 1$,  we have $D\leq 4^{\alpha }(1-t_{0})\|\varphi_{t}\|_{C^{0}}^{\alpha }$, then
\begin{eqnarray}
(1-t_{0})\geq \frac{D}{4^{\alpha }\|\varphi_{t_{0}}\|_{C^{0}}^{\alpha }}.
\end{eqnarray}
On the second case, we may assume that $1-t_{0}<\frac{A^{-1}}{2}$, the inequality implies that
\begin{eqnarray}
\|\varphi_{t_{0}}\|_{C^{0}}\leq 2\|\varphi_{1}\|_{C^{0}} +2 ,
\end{eqnarray}
then
\begin{eqnarray}
(1-t_{0})\geq \frac{D}{4^{\alpha }(2\|\varphi_{1}\|_{C^{0}}+2)^{\alpha }}.
\end{eqnarray}
Since $\sup \varphi_{1} \cdot \inf \varphi_{1}\leq 0$, we always have the following inequality
\begin{eqnarray}\label{PP1}
\begin{array}{lll}
(1-t_{0})&\geq & \frac{C'}{(\|\varphi_{1}\|_{C^{0}}+1)^{\alpha }}\\
&\geq & \frac{C'}{( Osc (\varphi_{1} )+1)^{\alpha }},\\
&= & \frac{C'}{( Osc (\phi )+1)^{\alpha }},\\
\end{array}
\end{eqnarray}
where $C'$ is a positive constant depending only on $(M, g_{SE})$. On the other hand, using proposition 5.4 again, we have
\begin{eqnarray}\label{PP2}
\begin{array}{lll}
(1-t_{0})\|\varphi_{1}-\varphi _{t_{0}}\|_{C^{0}} &\leq & (1-t_{0})^{2}A \|\varphi_{t_{0}}\|_{C^{0}}+1\\
&\leq & A D^{\frac{1}{\alpha }} +1 .\\
\end{array}
\end{eqnarray}
By inequalities (\ref{F3}), (\ref{PP1}) and (\ref{PP2}), we obtain
\begin{eqnarray}\label{MMT1}
F_{d\eta_{SE} }(\phi) \geq  \tilde{C}_{1} \frac{J_{d\eta_{SE} }(\phi)}{( Osc (\phi )+1 )^{\alpha }} -\tilde{C}_{2},
\end{eqnarray}
for all $\phi \in \mathcal H (\xi , \eta_{SE}, \Phi_{SE} , g_{SE})$, where $\tilde{C}_{1}$ and $\tilde{C}_{2}$ are positive constants depending only on the geometry of $(M, g_{SE})$.

Since $\varphi_{t}-\varphi_{1} \in \mathcal H (\xi , \eta_{SE}, \Phi_{SE} , g_{SE})$ and $\rho^{T}_{d\eta_{t}}\geq t(m+1)d\eta_{t}$, we can use corollary 2.8 to obtain the following estimate
\begin{eqnarray}\label{ZZ1}
Osc (\varphi_{t}-\varphi_{1} ) \leq I_{d\eta_{SE} }(\varphi_{t}-\varphi_{1} ) +\bar{C}(M, g_{SE}),
\end{eqnarray}
for $t\in [\frac{1}{2} , 1]$, where $\bar{C}(M, g_{SE})$ is a constant depending only on $(M, g_{SE})$.
By (\ref{MMT1}) and (\ref{ZZ1}), we have
\begin{eqnarray}\label{MMT2}
F_{d\eta_{SE} }(\varphi_{t}-\varphi_{1}) \geq  \tilde{C}_{3} \frac{J_{d\eta_{SE} }(\varphi_{t}-\varphi_{1})}{(J_{d\eta_{SE} }(\varphi_{t}-\varphi_{1})+1 )^{\alpha }} -\tilde{C}_{2},
\end{eqnarray}
for $t\in [\frac{1}{2} , 1]$, where $\tilde{C}_{3}$ is a constant depending only on $(M, g_{SE})$.

\medskip

 By the cocycle property of the functional $F$, formulas (\ref{K1}), (\ref{K2}), (\ref{C0}), nondecreasing of $(I_{d\eta }-J_{d\eta })(\varphi_{t})$ and the concavity of the log function, we have
\begin{eqnarray}\label{FF1}
\begin{array}{lll}
&& F_{d\eta_{SE} } (\varphi_{t}-\varphi_{1} )=F_{d\eta} (\varphi_{t})-F_{d\eta } (\varphi_{1} )\\
&=& \frac{-1}{t}\int_{0}^{t}(I_{d\eta }-J_{d\eta })(\varphi_{s}) ds  +\int_{0}^{1}(I_{d\eta }-J_{d\eta })(\varphi_{s}) ds \\ &&-\frac{1}{m+1}\log \{\frac{1}{V}\int_{M}e^{(t-1)(m+1)\varphi_{t} }(d\eta_{\varphi _{t}}
)^{m}\wedge \eta\}\\
&\leq & \frac{t-1}{t}\int_{0}^{t}(I_{d\eta }-J_{d\eta })(\varphi_{s}) ds  +\int_{t}^{1}(I_{d\eta }-J_{d\eta })(\varphi_{s}) ds\\
&& +\frac{(1-t)}{V}\int_{M}\varphi_{t}(d\eta_{\varphi_{t}}
)^{m}\wedge \eta \\
&=&\int_{t}^{1}(I_{d\eta }-J_{d\eta })(\varphi_{s}) ds-(1-t)(I_{d\eta }-J_{d\eta })(\varphi_{t})\\
&\leq & (1-t)\{(I_{d\eta }-J_{d\eta })(\varphi_{1})-(I_{d\eta }-J_{d\eta })(\varphi_{t})\}\\
&\leq & m(1-t) Osc (\varphi_{1}-\varphi_{t})\\
&\leq & m(1-t) \{I_{d\eta_{SE} }(\varphi_{t}-\varphi_{1}) +\frac{C_{1}(m)}{t}+C_{2}\}\\
&\leq & m(1-t) \{(m+1)J_{d\eta_{SE} }(\varphi_{t}-\varphi_{1}) +\frac{C_{1}(m)}{t}+C_{2}\}
\end{array}
\end{eqnarray}

By a same discussion in \cite{PSSW} (p1083), we know that (\ref{F3}), (\ref{ZZ1}), (\ref{MMT2}) and (\ref{FF1}) imply the Moser-Trudinger inequality (\ref{MT2}). We write out the proof in details just for reader's convenience.

Combining (\ref{MMT2}) with (\ref{FF1}), we have
\begin{eqnarray}\label{MMT3}
m(m+1)(1-t)J(t) +\tilde{C}_{4}(1-t) \geq  \tilde{C}_{3} \frac{J(t)}{(J(t)+1 )^{\alpha }} -\tilde{C}_{2},
\end{eqnarray}
for $t\in [\frac{1}{2}, 1]$, where $\tilde{C}_{4}$ is a constant depending only on $(M, g_{SE})$. Here we denote $J_{d\eta_{SE} }(\varphi_{t}-\varphi_{1})$ by $J(t)$ just for simplicity. (\ref{MMT3}) can also be written as
\begin{eqnarray}\label{Phong}
\frac{J(t)}{(J(t)+1)^{\alpha }}(\tilde{C}_{5}-(1-t)(J(t)+1)^{\alpha })\leq \tilde{C}_{6}(1-t)+\tilde{C}_{7}
\end{eqnarray}
where $\tilde{C}_{5}$, $\tilde{C}_{6}$ and $\tilde{C}_{7}$ are constants depending only on $(M, g_{SE})$.
We can suppose that there exists a $t' \in [\frac{1}{2} , 1]$ with
\begin{eqnarray}
(1-t')(J(t')+1)^{\alpha } =\frac{1}{2}\tilde{C}_{5}.
\end{eqnarray}
If not then we must have $(1-t)(J(t)+1)^{\alpha } <\frac{1}{2}\tilde{C}_{5}$ for all $t\in [\frac{1}{2} , 1]$. It would follow that
$J(\frac{1}{2})\leq \tilde{C}_{5}^{\frac{1}{\alpha }}$,  then (\ref{F3}) and (\ref{ZZ1}) imply (\ref{MT2}). Otherwise, from (\ref{Phong}) we have that $J(t')\leq \tilde{C}_{8}$ and $1-t' \geq \tilde{C}_{9}$, these also imply (\ref{MT2}).

\hfill $\Box$

Now, theorem 4.3 and theorem 6.1 imply the main theorem in the introduction.

\medskip

\section{A Miyaoka-Yau type inequality}
\setcounter{equation}{0}

{\bf Definition 7.1. } {\it Let  $(M, \xi , \eta , \Phi , g)$ be a compact Sasakian manifold with $[d\eta ]_{B} = \frac{2\pi }{m+1}
c_{1}^{B}(M, \mathcal F_{\xi})$. As above, we define $\mathcal S (\xi, \bar{J})$ to the space of all Sasakian structures which compatible with $(\xi , \eta , \Phi , g)$. Let's define two positive constants by
\begin{eqnarray*}
\alpha (\xi , \bar{J}) := \inf \{\lambda \quad | \quad 0\leq S^{T}_{d\eta '}\leq 2m\lambda \quad for \quad some \quad (\xi , \eta , \Phi , g)\in \mathcal S (\xi, \bar{J}) \};
\end{eqnarray*}
and
\begin{eqnarray*}
\beta (\xi , \bar{J}) := \sup \{\lambda \quad | \quad  S^{T}_{d\eta '}\geq 2m\lambda \quad for \quad some \quad  (\xi , \eta , \Phi , g)\in \mathcal S (\xi, \bar{J})\}.
\end{eqnarray*}}

\medskip

{\bf Remark: } {\it Since the mean value of transverse Scalar curvature $\bar{S}=2m(m+1)$ for any Sasakian structure in $\mathcal S (\xi, \bar{J})$, it is easy to see that $\alpha (\xi , \bar{J})\geq m+1$ and $0<\beta (\xi , \bar{J}) \leq m+1$. Obviously, if there exists a Sasakian-Einstein structure in $\mathcal S (\xi, \bar{J})$, then we have $\alpha (\xi , \bar{J})= m+1 = \beta (\xi , \bar{J})$.}

\medskip

{\bf Lemma 7.2. } {\it Let Let  $(M, \xi , \eta , \Phi , g)$ be a compact Sasakian manifold with $[d\eta ]_{B} = \frac{2\pi }{m+1}
c_{1}^{B}(M, \mathcal F_{\xi})$, and $(\xi , \eta', \Phi' , g') \in \mathcal S (\xi , \bar{J})$. Then we have
\begin{eqnarray}\label{MMY1}
\begin{array}{lll}
&&\int_{M} (2\pi )^{2} (2c_{2}^{B}(M, \mathcal F_{\xi })-\frac{m}{m+1}c_{1}^{B}(M, \mathcal F_{\xi })^{2})\wedge \frac{(\frac{1}{2}d\eta')^{m-2}}{(m-2)!}\wedge \eta' \\
&=& \int_{M} |Rm^{T}|^{2}-\frac{2(S^{T})^{2}}{m(m+1)}-\frac{(m-1)(m+2)}{m(m+1)}((S^{T})^{2}-(2m(m+1))^{2}) \frac{(\frac{1}{2}d\eta ')^{m}}{m!}\wedge \eta ',\\
\end{array}
\end{eqnarray}
where $Rm^{T}$ and $S^{T}$ are the transverse curvature tensor and the transverse scalar curvature of $(\xi , \eta', \Phi' , g')$.

}

\medskip

{\bf Proof. } By direct calculation, we have
\begin{eqnarray}
\begin{array}{lll}
&&\int_{M} (2\pi )^{2} (2c_{2}^{B}(M, \mathcal F_{\xi })-\frac{m}{m+1}c_{1}^{B}(M, \mathcal F_{\xi })^{2})\wedge \frac{(\frac{1}{2}d\eta')^{m-2}}{(m-2)!}\wedge \eta' \\
&=&\int_{M} \{ tr(Rm^{T}\wedge Rm^{T})-\frac{1}{m+1}tr Rm^{T}\wedge tr Rm^{T}\}\wedge \frac{(\frac{1}{2}d\eta')^{m-2}}{(m-2)!}\wedge \eta' \\
&=& \int_{M} |Rm^{T}|^{2}-|\rho^{T}|^{2}+\frac{1}{m+1}((S^{T})^{2}-|\rho^{T}|^{2}) \frac{(\frac{1}{2}d\eta')^{m}}{m!}\wedge \eta' \\
&=& \int_{M} |Rm^{T}|^{2}-(S^{T})^{2}-\frac{m+2}{m+1}(|\rho^{T}|^{2}-(S^{T})^{2}) \frac{(\frac{1}{2}d\eta')^{m}}{m!}\wedge \eta' .\\
\end{array}
\end{eqnarray}
On the other hand
\begin{eqnarray}
\begin{array}{lll}
&&\int_{M} (S^{T})^{2}-|\rho^{T}|^{2} \frac{(\frac{1}{2}d\eta')^{m}}{m!}\wedge \eta' \\
&=& \int_{M} \rho^{T}\wedge \rho^{T} \frac{(\frac{1}{2}d\eta')^{m-2}}{(m-2)!}\wedge \eta' \\
&=&\int_{M}4 m(m-1)(m+1)^{2} \frac{(\frac{1}{2}d\eta')^{m}}{m!}\wedge \eta' . \\
\end{array}
\end{eqnarray}
Combining the above two equalities, we get (\ref{MMY1}).

\hfill $\Box$

\medskip

In locally foliation chart $(x, z^{1}, \cdots , z^{m})$, setting
\begin{eqnarray}
Q_{i\bar{j}k\bar{l}}=R^{T}_{i\bar{j}k\bar{l}}-\frac{S^{T}}{m(m+1)}(g^{T}_{i\bar{j}}g^{T}_{k\bar{l}}+ g^{T}_{i\bar{l}}g^{T}_{k\bar{j}}).
\end{eqnarray}
It is easy to check that
\begin{eqnarray}\label{MMY2}
|Q|^{2}=|Rm^{T}|^{2}-\frac{2(S^{T})^{2}}{m(m+1)} .
\end{eqnarray}
 Combining (\ref{MMY1}) and (\ref{MMY2}), we have
\begin{eqnarray}\label{AAAA}
\begin{array}{lll}
&&\int_{M} (2\pi )^{2} (2c_{2}^{B}(M, \mathcal F_{\xi })-\frac{m}{m+1}c_{1}^{B}(M, \mathcal F_{\xi })^{2})\wedge \frac{(\frac{1}{2}d\eta')^{m-2}}{(m-2)!}\wedge \eta' \\
&\geq & \int_{M} -\frac{(m-1)(m+2)}{m(m+1)}((S^{T})^{2}-(2m(m+1))^{2}) \frac{(\frac{1}{2}d\eta')^{m}}{m!}\wedge \eta' . \\
\end{array}
\end{eqnarray}

Let's recall the Calabi functional on the space $\mathcal S (\xi , \bar{J})$, which was introduce by Boyer, Galicki and Simanca in \cite{BGS},
\begin{eqnarray}\label{DDDD}
\begin{array}{lll}
\mathcal Cal (\xi , \eta ' , \Phi', g')&=&\int_{M}(S_{d\eta' }^{T}-2m(m+1))^{2} (d\eta' )^{m}\wedge \eta'\\
&=&\int_{M}(S_{d\eta' }^{T})^{2}-(2m(m+1))^{2} (d\eta' )^{m}\wedge \eta' .\\
\end{array}
\end{eqnarray}

If $\inf_{\mathcal S (\xi , \bar{J})}\mathcal Cal =0$, for arbitrary $\epsilon >0$, we have a Sasakian structure $(\xi , \eta', \Phi' , g') \in \mathcal S (\xi , \bar{J})$ such that $\mathcal Cal (\xi , \eta ' , \Phi', g')\leq \epsilon $. Then, by (\ref{AAAA}), we have
\begin{eqnarray}\label{MMY3}
\begin{array}{lll}
&&\int_{M} (2\pi )^{2} (2c_{2}^{B}(M, \mathcal F_{\xi })-\frac{m}{m+1}c_{1}^{B}(M, \mathcal F_{\xi })^{2})\wedge \frac{(\frac{1}{2}d\eta')^{m-2}}{(m-2)!}\wedge \eta' \\
&\geq &  -\frac{(m-1)(m+2)}{m(m+1)}\epsilon .
\end{array}
\end{eqnarray}
Since $\epsilon $ is arbitrary, (\ref{MMY3}) implies the following theorem.

\medskip

{\bf Theorem 7.3. }  {\it Let  $(M, \xi , \eta , \Phi , g)$ be a compact Sasakian manifold with $[d\eta ]_{B} = \frac{2\pi }{m+1}
c_{1}^{B}(M, \mathcal F_{\xi})$. If $\inf_{\mathcal S (\xi , \bar{J})}\mathcal Cal =0$, then we have the following Miyaoka-Yau type inequality
\begin{eqnarray}\label{MMY4}
\int_{M}  (2c_{2}^{B}(M, \mathcal F_{\xi })-\frac{m}{m+1}c_{1}^{B}(M, \mathcal F_{\xi })^{2})\wedge (d\eta)^{m-2}\wedge \eta \geq 0 .
\end{eqnarray}
}

On the other hand, if $\alpha (\xi , \bar{J}) =m+1$, for arbitrary $\epsilon >0$, we have a Sasakian structure $(\xi , \eta', \Phi' , g') \in \mathcal S (\xi , \bar{J})$ such that $0\leq S^{T}\leq 2m(m+1+\epsilon )$. By (\ref{DDDD}),  we have
\begin{eqnarray}
\mathcal Cal (\xi , \eta ' , \Phi', g')\leq 2(2m)^{2}(m+1)\epsilon +(2m)^{2}\epsilon^{2}.
\end{eqnarray}
Then, we have the following corollary.

\medskip

{\bf Corollary 7.4. }  {\it Let  $(M, \xi , \eta , \Phi , g)$ be a compact Sasakian manifold with $[d\eta ]_{B} = \frac{2\pi }{m+1}
c_{1}^{B}(M, \mathcal F_{\xi})$. If $\alpha (\xi , \bar{J}) =m+1$, then $\inf_{\mathcal S (\xi , \bar{J})}\mathcal Cal =0$. In particulary, we also have the Miyaoka-Yau type inequality (\ref{MMY4}).

}

\medskip

As that in \cite{Bando}, we have the following proposition.

\medskip

{\bf Proposition 7.5. }  {\it Let  $(M, \xi , \eta , \Phi , g)$ be a compact Sasakian manifold with $[d\eta ]_{B} = \frac{2\pi }{m+1}
c_{1}^{B}(M, \mathcal F_{\xi})$. If the $\mathcal K$ energy functional $\mathcal V_{d\eta }$ is bounded below on the space $\mathcal H (\xi , \eta , \Phi , g)$, then, for arbitrary $\epsilon >0$, $M$ admits a Sasakian structure $(\xi , \eta', \Phi' , g')$ compatible with $(\xi , \eta , \Phi , g)$ such that $|S^{T}_{d\eta' }-2m(m+1)|\leq \epsilon $. In particularly, $\alpha (\xi , \bar{J})= m+1 = \beta (\xi , \bar{J})$.}

\medskip

{\bf Proof. } By proposition 4.4, there exists a smooth family
of solution $\{\varphi_{t}\}$ of (\ref{MA2}) for $t\in (0, 1 )$. Let $f(t)=(1-t)( I_{d\eta }-J_{d\eta })(\varphi_{t})$, by (\ref{MF5}), we have
\begin{eqnarray}
\frac{d}{dt } f(t)+(1-t)^{-1}f(t)=\frac{-1}{2(m+1)}\frac{d}{dt} \mathcal V_{d\eta }(\varphi_{t}).
\end{eqnarray}
Since $\mathcal V_{d\eta }$ is bounded below, the above equality implies that there exists a sequence $t_{i}\rightarrow 1$ such that
$f(t_{i})\rightarrow 0$ as $i\rightarrow +\infty $.
From  (\ref{C0}) and (\ref{2.17}), we have
\begin{eqnarray}
\begin{array}{lll}
\|h_{d\eta_{t}}\|_{C^{0}}&\leq & Osc (h_{d\eta_{t}})=(1-t) Osc (\varphi_{t})\\
&\leq & (1-t)((m+1)( I_{d\eta }-J_{d\eta })(\varphi_{t}) +\frac{C_{1}(m)}{t}+C_{2}).
\end{array}
\end{eqnarray}
So, there exists a sequence $t_{i}\rightarrow 1$ such that
$\|h_{d\eta_{t_{i}}}\|_{C^{0}}\rightarrow 0$ as $i\rightarrow +\infty $.
On the other hand, considering
\begin{eqnarray}
\rho_{d\eta_{t}}^{T}=t(m+1)d\eta_{t}+(m+1)(1-t)d\eta \geq t(m+1)d\eta_{t},
\end{eqnarray}
for arbitrary $\epsilon >0$, we get a Sasakian structure $(\xi , \tilde{\eta }, \tilde{\Phi } , \tilde{g})\in\mathcal S(\xi, \bar{J})$ such that
$S^{T}_{d\tilde{\eta } }-2m(m+1)\geq -\epsilon $ and $\|h_{d\tilde{\eta }}\|_{C^{0}}<\epsilon $.

Let's consider the Sasakian-Ricci flow (\ref{SR2}) with the initial data $d\tilde{\eta}_{0}=d\tilde{\eta }$. Since the initial $h_{d\tilde{\eta}_{0}}$ satisfies $\|h_{d\tilde{\eta }_{0}}\|_{C^{0}}<\epsilon $ and $\triangle_{0}h_{d\tilde{\eta }_{0}}\geq -2\epsilon$, by lemma 5.1, we have
\begin{eqnarray}
\|h_{d\tilde{\eta }_{s}}\|_{C^{0}}<4 e^{2(m+1)}\epsilon , \quad for \quad s\in [0, 2];
\end{eqnarray}
\begin{eqnarray}
\sup_{M}|dh_{d\tilde{\eta }_{s}}|_{s}^{2}<8 e^{4(m+1)}\epsilon^{2} , \quad for \quad s\in [1, 2];
\end{eqnarray}
and
\begin{eqnarray}
\triangle_{s}h_{d\tilde{\eta }_{s}}\geq -2 e^{2(m+1)}\epsilon , \quad for \quad s\in [0, 2].
\end{eqnarray}
From (\ref{SRR2}) and (\ref{SRR6}), setting $a=\frac{1}{4m}$, we have
\begin{eqnarray}
\begin{array}{lll}
&&(\frac{\partial }{\partial s }-\frac{1}{4}\triangle_{s})(|dh_{d\tilde{\eta }_{s}}|_{s}+\epsilon a (s-1)\triangle_{s}h_{d\tilde{\eta }_{s}})\\
&\leq & (m+1)(|dh_{d\tilde{\eta }_{s}}|_{s}+\epsilon a (s-1)\triangle_{s}h_{d\tilde{\eta }_{s}}) +\epsilon a \triangle_{s}h_{d\tilde{\eta }_{s}} \\&&-
(1+ \epsilon a (s-1))|\partial_{B}\bar{\partial }_{B} h_{d\tilde{\eta }_{s}}|_{s}^{2}\\
&\leq & (m+1)(|dh_{d\tilde{\eta }_{s}}|_{s}+\epsilon a (s-1)\triangle_{s}h_{d\tilde{\eta }_{s}}) + \triangle_{s}h_{d\tilde{\eta }_{s}}(\epsilon a
-\frac{1+\epsilon a (s-1)}{4m}\triangle_{s}h_{d\tilde{\eta }_{s}})\\
\end{array}
\end{eqnarray}
where we have used the Cauchy-Schwarz inequality $(\frac{1}{2}\triangle_{s} h)\leq m|\partial_{B}\bar{\partial }_{B} h|_{s}^{2}$. Equivalently, we have
\begin{eqnarray}\label{MMY5}
\begin{array}{lll}
&&(\frac{\partial }{\partial s }-\frac{1}{4}\triangle_{s})\{e^{1-s}(|dh_{d\tilde{\eta }_{s}}|_{s}+\epsilon a (s-1)\triangle_{s}h_{d\tilde{\eta }_{s}})\}\\
&\leq & e^{1-s} \triangle_{s}h_{d\tilde{\eta }_{s}}(\epsilon a
-\frac{1+\epsilon a (s-1)}{4m}\triangle_{s}h_{d\tilde{\eta }_{s}}).\\
\end{array}
\end{eqnarray}
Then, (\ref{MMY5}) implies that $e^{1-s}(|dh_{d\tilde{\eta }_{s}}|_{s}+\epsilon a (s-1)\triangle_{s}h_{d\tilde{\eta }_{s}})\leq 16e^{4(m+1)}\epsilon^{2}$ for $s\in [1, 2]$. Otherwise at the point of $[1, 2]\times M$ where it fails to hold for the first time $1<t_{0}\leq 2$, we have $e^{1-s}\epsilon a (s-1)\triangle_{s}h_{d\tilde{\eta }_{s}}\geq 8e^{4(m+1)}\epsilon^{2}$ and then  $\triangle_{s}h_{d\tilde{\eta }_{s}}\geq 32m e^{4(m+1)}\epsilon$. But, from (\ref{MMY5}), we have $\triangle_{s}h_{d\tilde{\eta }_{s}} \leq \epsilon $ at the point, which is a contradiction. So, we have
\begin{eqnarray}
\triangle_{s}h_{d\tilde{\eta }_{s}} \leq 64me^{4m+5}\epsilon \quad for \quad s=2,
\end{eqnarray}
and then
\begin{eqnarray}
|S^{T}_{d\tilde{\eta}_{s} }-2m(m+1)|\leq 32me^{4m+5}\epsilon \quad for \quad s=2.
\end{eqnarray}

 \hfill $\Box$

\medskip

{\bf Corollary 7.6.}  {\it Let  $(M, \xi , \eta , \Phi , g)$ be a compact Sasakian manifold with $[d\eta ]_{B} = \frac{2\pi }{m+1}
c_{1}^{B}(M, \mathcal F_{\xi})$. If the $\mathcal K$ energy functional $\mathcal V_{d\eta }$ is bounded below on the space $\mathcal H (\xi , \eta , \Phi , g)$, then we have the Miyaoka-Yau type inequality (\ref{MMY4}).
}

\medskip

As an application of theorem 4.3 and lemma 7.2, we have the following proposition.

\medskip

{\bf Proposition 7.7. } {\it Let  $(M, \xi , \eta , \Phi , g)$ be a compact Sasakian manifold with $[d\eta ]_{B} = \frac{2\pi }{m+1}
c_{1}^{B}(M, \mathcal F_{\xi})$ and
\begin{eqnarray}\label{CON721}
\int_{M} (2c_{2}^{B}(M, \mathcal F_{\xi })-\frac{m}{m+1}c_{1}^{B}(M, \mathcal F_{\xi })^{2})\wedge \frac{(\frac{1}{2}d\eta)^{m-2}}{(m-2)!}\wedge \eta =0.
\end{eqnarray}
 If $F_{d\eta }$ (or $\mathcal V_{d\eta}$) is proper in the space $\mathcal H(\xi , \eta , \Phi , g)$, then there must exists a Sasakian metric  $(\xi , \eta', \Phi' , g') \in \mathcal S (\xi , \bar{J})$ with constant curvature $1$. Furthermore, if $M$ is simply connected, then $(M, g')$ is isometric to a unit sphere.}

\medskip

{\bf Proof. } By theorem 4.3, there exists a Sasakian-Einstein $(\xi , \eta', \Phi' , g') \in \mathcal S (\xi , \bar{J})$. By lemma 7.2, formula (\ref{MMY2}) and the condition (\ref{CON721}), we have
\begin{eqnarray}
Q_{i\bar{j}k\bar{l}}=R^{T}_{i\bar{j}k\bar{l}}-2(g^{T}_{i\bar{j}}g^{T}_{k\bar{l}}+ g^{T}_{i\bar{l}}g^{T}_{k\bar{j}}),
\end{eqnarray}
i.e. $(\xi , \eta', \Phi' , g')$ is of constant transverse holomorphic bisectional curvature. On the other hand, using the relation (\ref{TRCR}) of the transverse curvature tensor and the Riemann curvature tensor (or see proposition 7.2 in \cite{TB}), it is not hard to see that the Riemannian manifold $(M, g')$ is of constant curvature $1$.
\hfill $\Box$

\medskip

{\bf Acknowledgements }

  The paper was written while the
 author was visiting McGill University. He would like to thank
ZheJiang University for the financial support and to thank McGill
University for the hospitality. The author would also like to
thank Prof. PengFei Guan and Xiangwen Zhang for their useful
discussion and help.

\hspace{1.4cm}

\bigskip

\end{document}